%% file: FrequencyDomain.tex
\documentclass[final]{siamonline250211}

\usepackage{subfig}
\usepackage{array}
\usepackage{lipsum}
\usepackage{amsfonts}
\usepackage{graphicx}
\usepackage{epstopdf}
\usepackage{algorithmic}

\ifpdf
  \DeclareGraphicsExtensions{.eps,.pdf,.png,.jpg}
\else
  \DeclareGraphicsExtensions{.eps}
\fi

\usepackage{enumitem}
\setlist[enumerate]{leftmargin=.5in}
\setlist[itemize]{leftmargin=.5in}

\newcommand{\Eb}{\mathop{\mathbb{E}}}

\newsiamremark{remark}{Remark}
\newsiamremark{hypothesis}{Hypothesis}
\crefname{hypothesis}{Hypothesis}{Hypotheses}
\newsiamthm{claim}{Claim}

\headers{Basis Choices for Frequency Domain Independence Tests}{J. Shi, W. Wang, W. Zhang, H. Bao, S. Chavez,J. Huang, Y. Wu, and K. Zhang}

\title{Basis Choices for Frequency Domain Statistical Independence Tests and Algorithms for 
Algebraic Relation Extraction \thanks{Submitted to the editors DATE.}}
\author{Juan Shi\thanks{Department of Mathematics, University of North Carolina at Chapel Hill, 
  Chapel Hill, NC (\email{juans@live.unc.edu}, \email{sergio@unc.edu}, \email{huang@email.unc.edu}). }
\and Wenbo Wang\thanks{Department of Biostatistics, University of North Carolina at Chapel Hill, 
    Chapel Hill, NC (\email{wenbo@live.unc.edu}).}
\and Wan Zhang\thanks{Center for Forecasting Science, Academy of Mathematics and System Science, Chinese Academy of Sciences, Beijing, China (\email{wanzhang@amss.ac.cn}).}
\and Han Bao\footnotemark[2]
\and Sergio Chavez\footnotemark[2]
\and Jingfang Huang\footnotemark[2] \thanks{Corresponding Author.}
\and Yichao Wu\thanks{Department of Mathematics, Statistics and Computer Science, University of Illinois at Chicago, Chicago, IL (\email{yichaowu@uic.edu}).}
\and Kai Zhang\thanks{Department of Statistics and Operations Research, University of North Carolina 
    at Chapel Hill, Chapel Hill, NC (\email{zhangk@email.unc.edu}).}}

\usepackage{amsopn}

\begin{document}

\maketitle

% REQUIRED
\begin{abstract}
In this paper, we explore how different selections of basis functions impact the efficacy of 
frequency domain techniques in statistical independence tests, and study different algorithms for extracting
low-dimensional algebraic relations from dependent data. We examine a range of complete orthonormal bases 
functions including the Legendre polynomials, Fourier series, Walsh functions, and standard and 
nonstandard Haar wavelet bases. We utilize fast transformation algorithms to efficiently transform 
physical domain data to frequency domain coefficients. The main focuses of this paper are the 
effectiveness of different basis selections in detecting data dependency using frequency domain data, e.g., whether varying basis choices significantly influence statistical power loss for small data with large noise; and on the stability of different optimization formulations for finding proper algebraic relations when data are dependent. We present numerical results to demonstrate the 
effectiveness of frequency domain-based statistical analysis methods and provide guidance 
for selecting the proper basis and algorithm to detect a particular type of relations. 
%The main contributions of this paper are thorough studies of basis choices in the frequency 
%domain independence tests with improved statistical power and efficient algorithms for finding
%underlying algebraic relations.

\end{abstract}

% REQUIRED

% REQUIRED

\input{sec01Intro-06192024}
\input{sec02Frequency0418-JH}

\input{sec03analysis0420-JH}

\input{sec04Numerical_0503_Juan}
\input{sec05Summary0503}

\bibliographystyle{siamplain}
\bibliography{bib/huang,bib/KZ_BIB_031223}

\end{document}

%% file: sec01Intro-06192024.tex
\section{Introduction}
Spectral or frequency domain analysis is a well-developed technique in science and engineering.  
Classical results include the Laplace transform, Fourier transform and Fourier series, 
orthogonal polynomials, and wavelet analysis. Instead of function values (physical domain data), 
the frequency domain technique studies the relations of expansion coefficients (frequency domain data)
which may show simpler structures when compared with their physical domain counterparts. 
For example, convolutions in the physical domain become simple products in the frequency 
domain, and physical domain differential equation models become algebraic relations for 
the expansion coefficients in the frequency domain. 

In this paper, we study how different choices of basis functions impact the numerical properties of 
the frequency domain techniques when generalized to statistical independence tests. 
For a given paired data set $\{(x_n, y_n)\}_{n=1}^N$ 
sampled from an underlying distribution, where both $x_n$ and $y_n$ may be scalar random variables, 
vectors, or tensors, the physical domain data set $\{(x_n, y_n)\}_{n=1}^N$ 
can be transformed into new frequency domain data which can be utilized to address one of the 
most fundamental questions in statistics: whether the variables $x$ and $y$ are dependent. 
Once the variables are discovered to be dependent, the frequency domain data can also be 
utilized to identify the structure of the dependence or relation.

The concepts of dependence and correlations of random data sets have been studied for over a
century. Classical parametric dependence tests include Pearson's correlation test, 
$\chi^2$-test, and t-test. When the dependence satisfies a linear relationship, linear regression 
is a well-developed technique with extensive applications in physical, biomedical, 
financial, and social sciences. More recent developments 
include~\cite{Balakrishnan2019hypothesis,berrett2020optimal,cao2020correlations,chatterjee2021new,
chatterjee2024survey,deb2020measuring,Deb2021,geenens2022hellinger,genest2019testing,genest2005locally,
gretton2007kernel,heller2016multivariate,heller2012consistent,heller2016consistent,jinmatteson18,
kojadinovic2009tests,lee2019testing,ma2019fisher,pfister2018kernel,sejdinovic2013,shi2020rate,szekely2007,
zhang2019bet,zheng2012generalized,Zhu2017}, with several notable contributions from computer 
scientists~\cite{kinney2014equitability,reshef2011detecting,reshef2015empirical,reshef2015equitability}, 
reflecting the general attention to this problem from the broader Data Science community.
Once the variables are determined dependent, it is a common practice to find any low-dimensional or compressible 
features in the data. A partial list of existing compressible features and their extraction strategies 
include the ``bump" function $y=f(x)=\Eb(Y|X = x)$ and bump hunting techniques~\cite{friedman1999bump},  
approximation with smooth functions in generalized additive models~\cite{hastie2017generalized}, 
projection pursuit regression~\cite{friedman1981projection}, optimization in variational 
inference~\cite{blei2017variational}, minimum average variance estimation 
(MAVE)~\cite{xia2002adaptive}, sliced inverse regression~\cite{li1991sliced}, and kernel density 
reduction and estimation~\cite{kerneldimension,rosenblatt1956remarks}.

The ideas of frequency domain statistical data analysis have been studied 
in the past, some recent results include the binary expansion testing 
(BET) \cite{zhang2019bet}, binary expansion adaptive symmetry test (BEAST) \cite{zhang2021beauty}, and binary expansion linear effect (BELIEF) \cite{brown2025belief}. 
Both these techniques transform the physical domain data 
nonlinearly to a set of frequency domain coefficients of the Walsh (Walsh-Hadamard) 
expansion of the probability density function. The independence test is then 
performed on the expansion coefficients to extract interesting statistical features. 
In this paper, in addition to the Walsh basis, we study different basis choices for 
frequency domain statistical independence tests, including the Legendre polynomials, 
Fourier series, and standard and nonstandard Haar wavelet bases. We introduce fast 
transformation algorithms and present the mathematically equivalent formulations for
different statistical signatures in the transformed frequency domain. When the data 
follow low-dimensional dependence relations, we study how to stably find the ``bump" 
function $y=f(x)=\Eb(Y|X = x)$ in the frequency domain. By studying a simplified yet 
representative class of separable relations in the form of $H(x) = G(y)$, we also 
present some numerical algorithms and address their stability issues for finding 
more general algebraic relations. The contributions of this paper include 
comparisons of different choices of basis functions for frequency domain statistical 
independence tests, fast transformation algorithms, and more stable bump hunting 
algorithms in the frequency domain. Compared with existing physical domain techniques, 
the frequency domain statistical data analysis framework may provide new strategies
for more efficient statistical testing and compressible feature extraction
with improved statistical power and numerical stability.

This paper presents the basic ideas of the frequency domain statistical data analysis 
from an {\it applied analysis perspective}, and is organized as follows. 
In \cref{sec:frequency}, we briefly discuss classical frequency domain analysis concepts 
and demonstrate how to ``sample" the expansion coefficients and higher order moments
efficiently using fast transformation algorithms. 
In \cref{sec:dependency}, we present the frequency domain statistical signatures 
of two independent scalar $X$ and $Y$ random variables and present two numerical 
estimators for statistical independence tests. In \cref{sec:bump}, we show how to 
stably find the ``bump" function $y=f(x)=\Eb(Y|X = x)$ using frequency domain data. 
We then present a canonical-correlation analysis (CCA) algorithm 
and its equivalent frequency domain singular value decomposition algorithm for 
identifying more general separable low-dimensional structures in the form of $H(x)=G(y)$.
In \cref{sec:num}, numerical results are presented to validate our analysis and 
to demonstrate the algorithms' performance and numerical stability. We finally 
summarize our results and discuss a few ongoing research efforts in \cref{sec:summary}. 
%including introducing well-developed
%algebra, applied harmonic analysis, and machine learning strategies to improve
%the statistical power of existing statistical tests, more rigorous and flexible
%definitions of ``compressible" and low-dimensional structures in the data, and 
%numerical techniques for very high dimensional data sets.

%% file: sec02Frequency0418-JH.tex
\section{Transform Physical Domain Data to Frequency Domain}
\label{sec:frequency}
Consider a single variable function $f(x)$, $-1 \leq x \leq 1$, the classical spectral/frequency domain analysis considers its expansion in the form of 
\begin{equation}
  f(x)=\sum_{i=0}^{\infty} c_{i} \phi_i(x)
\end{equation}
where $\{\phi_i(x) \}_{i=0}^{\infty}$ is a complete orthonormal basis set. Instead of studying
the physical domain function values $f(x)$ directly, the frequency domain technique studies the 
expansion coefficients $c_i$. The data set $\{c_i\}_{i=0}^{\infty}$ is referred to as the spectral or frequency domain data. Clearly, the mapping from $x$ to $\phi_i(x)$ is in general 
nonlinear. 

\subsection{Choice of Basis in One Dimension}
For a single-variable function, the choice of the basis set $\{\phi_i(x) \}$ is a well-studied 
topic in applied analysis. Each basis set is related to an inner product space. 
In this paper, we focus on four commonly used basis sets in one dimension: the Legendre 
polynomials, Fourier series, Walsh functions, and Haar wavelet basis. 

\vspace{0.05in}
{\noindent \textbf{Legendre Polynomials:}} When the inner product is defined as 
$\langle f, g\rangle = \int_{-1}^1 f(x) g(x) dx$, an orthogonal polynomial basis set $\{ P_k(x) \}_{k=0}^{\infty} $ 
can be constructed using the Gram-Schmidt process. The first four of these Legendre
polynomials are given by $P_0(x)=1$, $P_1(x)=x$, $P_2(x)=\frac12(3x^2-1)$, and 
$P_3(x)=\frac12 (5x^3-3x)$. Similar to the properties of other orthogonal polynomial basis sets, the Legendre polynomials have a generating function
\begin{displaymath}
	\frac{1}{\sqrt{1-2xt+t^2}} = \sum_{n=0}^{\infty} P_n(x) t^n
\end{displaymath}
and a three-term recurrence relation  
\begin{displaymath}
  (n+1) P_{n+1}(x) = (2n+1) x P_n(x) - n P_{n-1}(x),
\end{displaymath}
which is often applied to accelerate the evaluation of the polynomial expansions. 
Legendre polynomials are solutions to Legendre's differential equation
$$(1-x^2) P''_n(x) - 2x P_n'(x) + n(n+1) P_n(x) =0.$$
Interested readers are referred to \cite{abramowitz1988handbook,lozier2003nist}
for a complete list of the properties of the Legendre (and other orthogonal) polynomials.
In this paper, we use orthonormal basis sets for all expansions. 
Instead of the standard Legendre polynomials which satisfy 
\begin{displaymath}
	\int_{-1}^{1} P_M(x) P_n(x) dx = \frac{2}{2n+1} \delta_{nm},
\end{displaymath}
the normalized orthonormal Legendre polynomial basis is given by 
$L_n(x)= \sqrt{\frac{2n+1}{2}} P_n(x).$

\vspace{0.05in}
{\noindent \textbf{Fourier Series:}} Fourier series is a useful tool to study 
smooth periodic functions. Using the complex representation, the Fourier series
of a function $f(x)$ which is periodic with a fundamental period of $2$ is
$$f(x) =\sum_{k=-\infty}^{\infty} c_k e^{ i k \pi x} $$
where $i=\sqrt{-1}$ and 
$c_k= \frac12 \int_{-1}^{1} f(x) e^{- i k \pi x} dx.$
Compared with orthogonal polynomials, the properties of the Fourier series are much simpler:
The generating function of the Fourier series is the (periodic) Dirac Delta function, 
the three-term recurrence relation is simply $e^{i(j+1) \pi x} = e^{ij \pi x} e^{i \pi x}$,
and the basis $e^{ij\pi x}$ is the solution of the simple differential equation
$$u''_n(x) + (j \pi )^2 u(x) =0.$$
In this paper, similar to the Legendre basis $L_n(x)$, we use the normalized complex form 
Fourier basis functions $\phi_n(x) = \frac{1}{\sqrt{2}} e^{i n \pi x}$. 

\vspace{0.05in}
{\noindent \textbf{Walsh Functions:}} Wavelet transforms are commonly used in
image and signal analysis. One particular wavelet basis set is the Walsh 
functions~\cite{beauchamp1976walsh,walsh1923closed} which consist of trains of 
square pulses (with the allowed states being $-1$ and $1$) such that transitions 
may only occur at fixed intervals of a unit time step. The initial state is always 
$+1$, and the functions satisfy the orthogonality relations and they all have norm
one. In particular, the $2^n$ Walsh functions of order $n$ are given by the rows 
of the Hadamard matrix $H_{2^n}$ when arranged in so-called ``sequency" 
order~\cite{thompson2017interferometry}. 
% There are $2^n$ Walsh functions of length $2^n$.
% illustrated in Fig.~\ref{fig:walsh} for $n=1$, $2$ and $3$. 
The Walsh functions are often considered as the discrete 
version of the $sine$ and $cosine$ functions in a Fourier series expansion.
The wavelet basis may perform better when the underlying function has discontinuities 
as discussed in \cite{strang1993wavelet}, where the author compared the performance 
of wavelet expansions with Fourier series for different applications. 
%\begin{figure}[htbp]
%\centering
%\includegraphics[width=14cm,height=10cm]{Figures/walsh_example.jpg}
% \caption{Illustration of Walsh functions for levels $1$ ($W(1,:)$), $2$ ($W(2,:)$) and $3$ ($W(3,:)$). 
% \label{fig:walsh}}
%\end{figure}

\vspace{0.05in}
{\noindent \textbf{Haar Wavelets:}} Another choice of wavelet basis is the Haar 
wavelet \cite{haar1910theorie} consisting of a sequence of square-shaped functions 
constructed using the scaling function $s(t)$ and mother wavelet function 
$m(t)$ given by 
$$s(x)= \left\{ \begin{array}{ll} 1 \quad & 0 \leq x <1 \\ 0 \quad & otherwise \end{array} \right. $$
$$m(x) = \left\{ \begin{array}{ll} 1 \quad & 0 \leq x < \frac12  \\ -1 \quad & \frac12 \leq x < 1 \\
	0 \quad & otherwise. \end{array} \right. $$
The non-normalized basis functions at level $i$ are given by the formula 
$ \psi_{i,j}(t) = m(2^{i}t - j)$ where $j = 0, 1, ..., 2^{i} - 1$.
%In Fig.~\ref{fig:haarbasis}, we plot the first eight non-normalized Haar basis functions (up to level $2$).
%The first two plots show the scaling function and mother wavelet function, the next two functions are the 
%level $1$ basis functions, and the last four functions are the level $2$ basis functions. 
%\begin{figure}[h!]
%    \centering
%    \includegraphics[width=14cm]{./Figures/haarbasis1.jpg}
%    \caption{The first 8 Haar basis functions}
%    \label{fig:haarbasis}
%\end{figure}
Similar to the Legendre and Fourier basis functions, we also normalize the Haar basis
functions using $H_0(t)=s(t)$ and $H_k(t)= 2^i \psi_{i,j}(t)$ so that the rescaled 
Haar basis functions satisfy $||H_k(x)||_2=1$ for all $k$ values.
An interesting feature of the Haar basis is the ``local" (sparse) property of the higher-level 
basis functions, i.e., the basis is only nonzero in a small region for large $k$ values.

\subsection{Choice of Basis in Higher Dimensions}
The single variable frequency domain analysis can be easily extended to study multi-variable functions. 
Using the tensor/Kronecker product, higher-dimensional basis sets can be constructed using one-dimensional basis functions. For a function $f(x,y)$ with two variables, the basis functions are $\{ \phi_i(x) \psi_j(y) \}_{i,j=0}^{\infty}$ where $\{ \phi_i(x) \}$ and 
$ \{ \psi_j(y) \}$ are the (not necessarily the same) one-dimensional basis 
sets for the $x$- and $y$-variables, respectively. The function expansion is 
then given by 
$$ f(x,y)=\sum_{i,j} c_{i,j} \phi_i(x) \psi_j(y). $$
Note that non-tensor-product bases can also be used and a more general expansion is then given by 
$$ f(x,y)=\sum_{k} c_{k} \Phi_k(x,y)$$
where $\{ \Phi_k(x,y) \}_{k=0}^{\infty}$ is an arbitrary two-dimensional orthonormal 
basis. 

\begin{figure}[htbp]
    \centering
    \includegraphics[width=10cm]{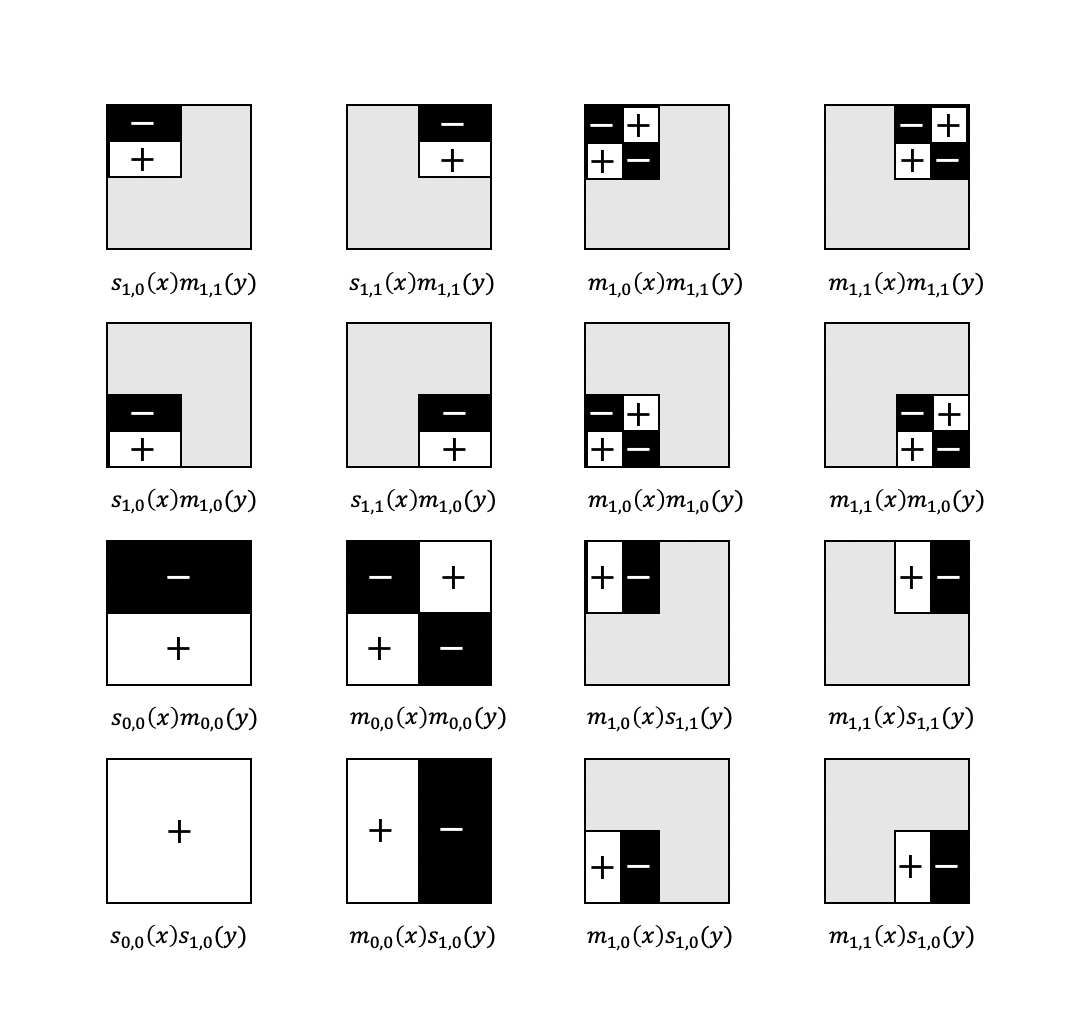}    \caption{The nonstandard construction of a two-dimensional Haar wavelet basis}
 \label{fig:nonstandard_haarbasis}
\end{figure}

An interesting non-tensor-product wavelet basis is the nonstandard Haar wavelet 
basis \cite{stollnitz1995wavelets} constructed using the one-dimensional 
functions $s_{i,j}(t) = s(2^{i}t - j)$ and $m_{i,j}(t) = m(2^{i}t - j)$ with sparser and more local properties. 
In Fig.~\ref{fig:nonstandard_haarbasis}, we plot the first 2 levels of 
nonstandard Haar basis functions in two dimensions. In the plot, the ``+" and ``$-$" 
signs represent a positive constant $+c$ and a negative constant $-c$, respectively. 
The constant $c$ is chosen so that the nonstandard Haar wavelet 
basis functions are normalized. At deeper levels, the basis functions are zero in 
most regions and are hence ``sparse" and ``local".
Note that the nonstandard Haar basis is naturally associated with an adaptive hierarchical tree structure, which allows faster computation of the expansion coefficients. For high-dimensional data located on a low-dimensional manifold, the ``sparse" nonstandard Haar wavelet basis is a better choice compared 
with other ``all-grid" tensor-product based basis choices such as the 
Legendre polynomials, Fourier series, and Walsh functions.

\subsection{Transform from Physical Domain to Frequency Domain}
\label{sec:transformed}
We demonstrate the transformation from the physical domain to the frequency domain using two-dimensional data sets. Assume that $N$ sampled data points 
$\{ (x_n,y_n) \}_{n=1}^N$ are located in the square $[a,b]^2$,
where $[a,b]=[-1,1]$ for the Legendre polynomials and Fourier series, 
and $[a,b]=[0,1]$ for the Walsh and Haar wavelet bases.
We consider the two-dimensional basis set 
$\{ \phi_{k_1}(x) \phi_{k_2}(y) \}_{k_1, k_2 = 0}^{\infty}$ constructed using the 
tensor product applied to a general one-dimensional orthonormal basis set 
$\{ \phi_k(x) \}_{k=0}^{\infty}$ for the domain $[a,b]$. We use the same basis set for both $x$- and $y$-variables. The generalizations of our results to a two-dimensional basis set constructed using the tensor product of two different 
one-dimensional basis sets and to a non-tensor-product orthonormal basis are 
straightforward. We skip the details in this paper.

Assume that the frequency domain representation of the underlying physical domain
probability density function $f(x,y)$ for the data points $\{ (x_n,y_n) \}_{n=1}^N$ 
is given by 
$$ f(x,y)=\sum_{k_1,k_2} c_{k_1,k_2} \phi_{k_1} (x) \phi_{k_2} (y).$$
For fixed $k_1$ and $k_2$, by applying the orthogonality properties of the basis functions, 
we have the following theorem which states that each datum $\phi_{k_1}(x_n) \phi_{k_2}(y_n)$ 
can be considered as a ``sample" of the expansion coefficient $ c_{k_1,k_2}$. 
\begin{theorem}
%\begin{equation}			
	$\Eb [\phi_{k_1}(X) \phi_{k_2}(Y)] = \int  \int \phi_{k_1}(x) \phi_{k_2}(y) f(x,y) dx dy = c_{k_1,k_2}.$
%\end{equation}
\end{theorem}
We introduce the following definition.
\begin{definition}
	For each pair of physical domain data $(x_n, y_n)$, the (nonlinearly) transformed data sets  
	$\{ \phi_{k_1}(x_n) \phi_{k_2}(y_n) \}_{n=1}^N$ for $k_1, k_2 = 0, \cdots, \infty$ are called the frequency domain data sets.
\end{definition}
%Note that the transformation from the physical domain to the frequency domain utilizes the nonlinear basis 
%functions $\phi_{k_1}(x) \phi_{k_2}(y)$ when $k_1, k_2 \geq 2$.
Given $N$ sample points, using the frequency domain data sets $\{ \phi_{k_1}(x_n) \phi_{k_2}(y_n) \}_{n=1}^N$
and Monte Carlo integration, the numerical estimates of the frequency domain coefficients are 
\begin{equation} 
	\hat{c}_{k_1,k_2} = \frac{1}{N} \sum_{n=1}^N \phi_{k_1}(x_n) \phi_{k_2}(y_n).
\end{equation}
The frequency domain coefficient matrix $C=[c_{k_1,k_2}]$ can be approximated by $\hat{C}=[\hat{c}_{k_1,k_2}]$. 

When a tensor product is applied to construct a two-dimensional orthonormal basis, we can also 
study the statistical properties of the special sets of frequency domain data 
$\{ \phi_{k_1}(x_n) \}_{n=1}^N $ and  $\{ \phi_{k_2}(y_n) \}_{n=1}^N.$
Clearly, 
$$\Eb [ \left[ \begin{array}{c} \phi_0(X) \\ \phi_1(X) \\ \vdots \\ \phi_K(X) \end{array} \right]] = 
	  \left[ \begin{array}{c} c_{0,0} \\ c_{1,0} \\ \vdots \\ c_{K,0} \end{array} \right] 
		  \mbox{ and }  
   \Eb [ \left[ \begin{array}{c} \phi_0(Y) \\ \phi_1(Y) \\ \vdots \\ \phi_K(Y) \end{array} \right] ] =
          \left[ \begin{array}{c} c_{0,0} \\ c_{0,1} \\ \vdots \\ c_{0,K} \end{array} \right].$$
Following classical statistical analysis theory, we define the variance matrices of 
$\{ \phi_{k_1}(x_n) \}_{n=1}^N $ and  $\{ \phi_{k_2}(y_n) \}_{n=1}^N$ as 
\begin{eqnarray*}
	\Sigma_{XX} & =  & \Eb [\left( 
	  \left[ \begin{array}{c} \phi_0(X) \\ \phi_1(X) \\ \vdots \\ \phi_K(X) \end{array} \right] - 
		  \left[ \begin{array}{c} c_{0,0} \\ c_{1,0} \\ \vdots \\ c_{K,0} \end{array} \right] \right) 
	\left( 
	  \left[ \begin{array}{c} \phi_0(X) \\ \phi_1(X) \\ \vdots \\ \phi_K(X) \end{array} \right] - 
		  \left[ \begin{array}{c} c_{0,0} \\ c_{1,0} \\ \vdots \\ c_{K,0} \end{array} \right] \right)^T
			  ] \\
	& = & \Eb [ \phi_k(X) \phi_m(X)] - \left[ \begin{array}{c} c_{0,0} \\ c_{1,0} \\ \vdots \\ c_{K,0} \end{array} \right] 
		[ c_{0,0}, c_{1,0}, \cdots, c_{K,0} ], \\
	\Sigma_{YY} & =  & \Eb [\left( 
	  \left[ \begin{array}{c} \phi_0(Y) \\ \phi_1(Y) \\ \vdots \\ \phi_K(Y) \end{array} \right] - 
		  \left[ \begin{array}{c} c_{0,0} \\ c_{0,1} \\ \vdots \\ c_{0,K} \end{array} \right] \right) 
	\left( 
	  \left[ \begin{array}{c} \phi_0(Y) \\ \phi_1(Y) \\ \vdots \\ \phi_K(Y) \end{array} \right] - 
		  \left[ \begin{array}{c} c_{0,0} \\ c_{0,1} \\ \vdots \\ c_{0,K} \end{array} \right] \right)^T
			  ] \\
	& = & \Eb [ \phi_k(Y) \phi_m(Y)] - \left[ \begin{array}{c} c_{0,0} \\ c_{0,1} \\ \vdots \\ c_{0,K} \end{array} \right] 
				  [ c_{0,0}, c_{0,1}, \cdots, c_{0,K} ]. 
\end{eqnarray*}
As will be explained in Sec.~\ref{sec:alg}, we can study the relations between the physical domain
$x$- and $y$-variables using the frequency domain covariance matrix of 
$\{ \phi_{k_1}(x_n) \}_{n=1}^N $ and  $\{ \phi_{k_2}(y_n) \}_{n=1}^N$ defined as
\begin{eqnarray*}
	\Sigma_{XY} & =  & \Eb [\left( 
	  \left[ \begin{array}{c} \phi_0(X) \\ \phi_1(X) \\ \vdots \\ \phi_K(X) \end{array} \right] - 
		  \left[ \begin{array}{c} c_{0,0} \\ c_{1,0} \\ \vdots \\ c_{K,0} \end{array} \right] \right) 
	\left( 
	  \left[ \begin{array}{c} \phi_0(Y) \\ \phi_1(Y) \\ \vdots \\ \phi_K(Y) \end{array} \right] - 
		  \left[ \begin{array}{c} c_{0,0} \\ c_{0,1} \\ \vdots \\ c_{0,K} \end{array} \right] \right)^T
			  ] \\
	& = & \Eb [ \phi_k(X) \phi_m(Y)] - \left[ \begin{array}{c} c_{0,0} \\ c_{1,0} \\ \vdots \\ c_{K,0} \end{array} \right]
		[ c_{0,0}, c_{0,1}, \cdots, c_{0,K} ]. 
\end{eqnarray*}
Note that $\Sigma_{XY}$ is the frequency domain coefficient matrix $C=[c_{k_1,k_2}]$ minus a rank one
perturbation matrix.  As we always choose $\phi_0(x)$ to be the constant function and perform the copula transformation, the first column and 
first row of the matrices $\Sigma_{XX}$, $\Sigma_{YY}$, and $\Sigma_{XY}$ are therefore always zero except for the first entry, 
\begin{displaymath}
	\Sigma_{XX} =  \left[ \begin{array}{cccc} 1 & 0 & \cdots & 0 \\
		0 & & & \\
		\vdots & & \tilde{\Sigma}_{XX} & \\
		0 & & & \\
	\end{array} \right], \quad 
	\Sigma_{YY}  = \left[ \begin{array}{cccc} 1 & 0 & \cdots & 0 \\
		0 & & & \\
		\vdots & & \tilde{\Sigma}_{YY} & \\
		0 & & & \\
	\end{array} \right], \quad
	\Sigma_{XY}  =  \left[ \begin{array}{cccc} 1 & 0 & \cdots & 0 \\
		0 & & & \\
		\vdots & & \tilde{\Sigma}_{XY} & \\
		0 & & & \\
	\end{array} \right]. \\
\end{displaymath}
We therefore only need to focus on $\tilde{\Sigma}_{XX}$, $\tilde{\Sigma}_{YY}$, and $\tilde{\Sigma}_{XY}$.
Finally, we define the variance of the frequency domain data $\phi_{k_1}(x) \phi_{k_2}(y)$ as
$$ var[\phi_{k_1}(X) \phi_{k_2}(Y)]  =  \Eb [((\phi_{k_1}(X) \phi_{k_2}(Y) - c_{k_1,k_2})^2] = 
  \Eb [((\phi_{k_1}(X) \phi_{k_2}(Y))^2] -  c_{k_1,k_2}^2. $$
Higher order moments of the frequency domain data can also be computed which may carry interesting
information on the statistical properties of the physical domain data.

\vspace{0.05in}
{\noindent \bf Comment:} For non-tensor-product bases, e.g., the nonstandard Haar wavelet basis
functions, the matrices $\tilde{\Sigma}_{XX}$, $\tilde{\Sigma}_{YY}$, and $\tilde{\Sigma}_{XY}$
are unfortunately not defined. However, it is still possible to define the expectations and
higher order moments of the frequency domain data $\{ \Phi_k(x_n,y_n) \}$ for $n=1,\cdots,N$ and $k=0, \cdots, \infty$.
%It becomes challenging to study the relations between
%$x$- and $y$-variables using the frequency domain covariance matrix $\tilde{\Sigma}_{XY}$.

One of the main contributions of this paper is a more systematic study of the statistical properties of 
the physical domain data $(x_n, y_n)$  using the transformed frequency domain data including 
the (numerically computed random) {\bf frequency domain coefficient matrix} $\hat{C}_{(K+1) \times (K+1)}$,
{\bf variance matrices} $\tilde{\Sigma}_{XX}$, $\tilde{\Sigma}_{YY}$, 
$\tilde{\Sigma}_{XY}$, $var[\phi_{k_1}(X) \phi_{k_2}(Y)]$, and higher order moments. Our current study focuses on the low-frequency domain information, where the frequency domain
expansion is truncated after the first $K$ terms with $K \leq c \cdot \sqrt{N}$, $c<1$. 
Intuitively, accurately capturing all the high-frequency coefficients $\hat{c}_{k_1,k_2}$ 
would require enough sampling points in one wavelength for each variable. 
However, it is still possible to get certain high-frequency information when the 
information is sufficient in certain perspectives (e.g., physical data are dense 
enough in a small spatial region). This paper studies how to extract statistical information 
from the physical and frequency domain data under different scenarios.

\subsection{Fast Transformation Algorithms}
The coefficient and variance matrices $\hat{C}$, $\tilde{\Sigma}_{XX}$, $\tilde{\Sigma}_{YY}$, 
$\tilde{\Sigma}_{XY}$, and $var[\phi_{k_1}(X) \phi_{k_2}(Y)]$ can be efficiently computed 
using fast transformation algorithms. Using the Fourier basis as an example, we demonstrate how to efficiently compute all the entries in the coefficient matrix $\hat{C}$ 
$$\hat{c}_{k_1,k_2} = \frac{1}{N} \sum_{n=1}^N e^{i k_1 \pi x_n}  e^{ i k_2 \pi y_n}  $$
for $k_1, k_2 = -K, \cdots, K$. As the data points $(x_n, y_n)$ are in general not uniformly distributed, and evaluating the summation directly would require $O(N K^2)$ operations. 
In two dimensions, when $K \sim O(\sqrt{N})$, the total amount of work in a direct method is 
approximately $O(N^2)$. Following the notation in \cite{greengard2004accelerating}, we notice 
that the summation for computing 
$\hat{c}_{k_1,k_2}$ is a special case of the Type I nonuniform discrete Fourier transform (see Eq.~(1) in 
\cite{greengard2004accelerating}) where the weight for each data pair is $1$, therefore existing nonuniform fast Fourier transform (NuFFT) algorithms can be applied to reduce the total amount of work 
from $O(N^2)$ to $O(N \log N)$~\cite{debroy2008accelerating,
dutt1993fast,greengard2004accelerating,lee2005type,potts2001fast}.
The NuFFT algorithm is based on the $O(N \log N)$ Fast Fourier Transform (FFT) 
algorithm for uniformly distributed points~\cite{cooley1965algorithm}. Interested readers 
are referred to \cite{heideman1984gauss} for the history of the development of FFT algorithms. 

For other basis choices, existing uniform fast Legendre transform (FLT) packages are at least 
$5 \sim 10$ times slower than the well-developed (and hardware-implemented) uniform FFT solvers; both uniform and non-uniform fast Walsh-Hadamard transform (FWHT) and fast Haar 
transform (FHT) are available to accelerate the computation of the random frequency domain 
coefficient matrix $\hat{C}$. When the nonstandard Haar wavelet basis is used, because of the sparsity
and locality in the basis functions and by utilizing the associated hierarchical tree structure, the accelerated computation of the nonstandard Haar wavelet expansion can be asymptotically 
optimal in efficiency and allows for easy coupling with an adaptive tree structure when data 
are located on low-dimensional manifolds. This feature makes the nonstandard
Haar basis appealing for very high-dimensional data sets which are located on a low-dimensional
manifold. In engineering applications, the Walsh and Haar wavelet bases may show better 
performance when the underlying function has discontinuities as discussed in \cite{strang1993wavelet}.
We have implemented and tested some of these fast transformation algorithms. Our numerical results 
show that the fast transformations from the physical domain to the frequency domain can be orders of magnitude more efficient, especially for a large number of data points.

%% file: sec03analysis0420-JH.tex
\section{Frequency Domain Statistical Analysis}
\label{sec:alg}
We demonstrate how the frequency domain data can be utilized 
to study the dependence of two scalar random variables $X$ and $Y$, and when they 
are dependent, how to find possible low-dimensional structures in the data. 

\subsection{Independence Test}
\label{sec:dependency}
In the physical domain, two random variables $X$ and $Y$ are independent if
and only if their joint distribution function satisfies $f(x,y)=h(x) g(y)$, where
$h(x)$ and $g(y)$ are the marginal distribution functions for $X$ and 
$Y$, respectively. When a two-variable function $f(x,y)$ is given by
the product of two single variable functions $h(x)$ and $g(y)$, we refer
to $f(x,y)$ as a {\bf rank-1} function. 

In the frequency domain, when using the same set of one-dimensional basis functions 
$\{ \phi_k \}_{k=0}^{\infty}$ with $\phi_0=\frac{1}{||1||}$ (constant function), 
one can expand $h(x)$ and $g(y)$ as
$$h(x)=\sum_{j=0}^{\infty} h_j \phi_j(x) \mbox{ and } 
g(y)=\sum_{k=0}^{\infty} g_k \phi_k(y),$$
respectively. The frequency domain expansion of $f(x,y)$ is then given by
$$ f(x,y)= \sum_{j,k=0}^{\infty} c_{j,k}  \phi_j(x) \phi_k(y)=\sum_{j,k=0}^{\infty} h_j g_k  \phi_j(x) \phi_k(y). $$  
Therefore, the coefficients satisfy the relation $c_{j,k}= h_j \cdot g_k$ and 
the truncated coefficient matrix 
$$C_{K \times K} = [h_0, h_1, \cdots, h_K]^T [g_0, g_1, \cdots, g_K] $$
is a rank-1 matrix. We have the following theorem for the frequency domain
data.
\begin{theorem}
	\label{theorem:ind01}
Assume the joint probability distribution function $f(x,y)$ of two random variables 
$X$ and $Y$ has the frequency domain expansion 
	$f(x,y)=\sum_{j,k=0}^{\infty} c_{j,k}  \phi_j(x) \phi_k(y)$. 
Then $X$ and $Y$ are independent if and only if the expansion coefficient matrix $C=[c_{j,k}]$ 
is a rank-1 matrix $C = [h_0, h_1, \cdots]^T [g_0, g_1, \cdots]$.
\end{theorem}

Determining whether a random matrix is rank-1 is theoretically possible but numerically more challenging. 
Therefore, following standard practice in statistical analysis, to study whether $Y$ is dependent 
on $X$, one can first apply a ``copula" transform to $Y$ so that the transformed $Y$ variable
follows a uniform distribution on the interval $[0,1]$. We can simplify Theorem~\ref{theorem:ind01}
as follows.
\begin{theorem}
\label{theorem:ind02}
Assume the joint probability distribution function $f(x,y)$ of two random variables 
$X$ and $Y$ has the frequency domain expansion $f(x,y)=\sum_{j,k=0}^{\infty} c_{j,k}  \phi_j(x) \phi_k(y)$
and the marginal distribution of $Y$ follows the uniform distribution $U[0,1]$.
Then $X$ and $Y$ are independent if and only if the expansion coefficient matrix $C=[c_{j,k}]$ 
is a rank-1 matrix $C = [h_0, h_1, h_2, \cdots]^T [1, 0, 0, \cdots]$. Equivalently,
the expansion coefficients $c_{j,k}=0$ for any $k \neq 0$.
\end{theorem}
Therefore, after applying the copula transform to $y$-variable, the task of testing independence
becomes that of testing whether the frequency domain data satisfy $c_{j,k}=0$ for any $k \neq 0$.

When $Y$ and $X$ are independent and after performing the copula transform to $Y$, we have 
$\Eb [\phi_j(X)\phi_k(Y)]=c_{j,k}=0$ ($k \neq 0$) independent of the marginal distribution of $X$. However, 
other statistical properties of $\phi_j(x)\phi_k(y)$ will depend on the marginal 
distribution of $X$ (e.g., its variance). This dependency implies that for any frequency domain 
estimator in the statistical independence test, its distribution for the independent case cannot be 
precomputed as the marginal distribution of $X$ is unknown. It is possible to run a permutation
test which permutes either the $X$, or  $Y$, or both variables 
\cite{annis2005permutation,collingridge2013primer,edgington2007randomization,fisher1936design}, 
and the permuted data sets are used to estimate the distribution function of $\phi_j(x)\phi_k(y)$ for
the independent cases. The on-the-fly permutation test is normally very expensive. 
For optimal numerical stability and efficiency (and simplified notation 
and discussion), when studying the dependence between two random variables $X$ and $Y$, 
we assume that the cumulative distribution function of the physical domain data 
$\{ (x_n,y_n) \}_{n=1}^N$ is a copula and the marginal probability distributions 
of both variables $X$ and $Y$ are uniform on the interval $[0,1]$. 
Introducing a proper orthonormal basis set $\{ \phi_k(x) \}_{k=0}^{\infty}$ for the 
domain $[0,1]$ with  $\phi_0(x)=\frac{1}{||1||}$ and constructing the corresponding orthonormal basis for 
$[0,1]^2$ using the tensor product, Theorem~\ref{theorem:ind02} can then be simplified to 
the following theorem which provides a set of mathematically equivalent frequency domain signatures 
when $X$ and $Y$ are independent.
\begin{theorem}
	\label{theorem:ind03}
	Assume $X, Y \sim U[0,1]$. Then the following statements are equivalent. 
	\begin{enumerate}
		\item $X$ and $Y$ are independent. 
		\item The probability density function is the constant function $f(x,y)=1$, $0 \leq x, y \leq 1$.
		\item $c_{k_1,k_2}=0$ when $(k_1, k_2) \neq (0,0)$, where $c_{k_1,k_2}$ are the
			expansion coefficients of 
			$$f(x,y)=\sum_{k_1,k_2} c_{k_1,k_2} \phi_{k_1} (x) \phi_{k_2} (y),$$
			and $\{ \phi_k \}$ is a complete set of orthonormal basis with 
			$\phi_0=\frac{1}{||1||}.$
		\item For any $\vec{\alpha}$ and $\vec{\beta}$, define $\psi_1(x) = \sum_k \alpha_k \phi_k(x)$ 
			and $\psi_2(y) = \sum_k \beta_k \phi_k(y)$. Then, the expectation $\Eb (\psi_1(X) \psi_2(Y))=\alpha_0 \cdot \beta_0$.
		\item Corr$(\psi_1(X), \psi_2(Y)) = 0$ for any  $\vec{\alpha}$ and $\vec{\beta}$. 
	\end{enumerate}
\end{theorem}
We also list some necessary (but not sufficient) conditions for $X$ and $Y$ to be independent.
\begin{theorem}
	Assume $X, Y \sim U[0,1]$, $X$ and $Y$ are independent, and $\phi_0(x)=\frac{1}{||1||}$, 
	then the frequency domain data have the following features.
	\begin{enumerate}
		\item $\Eb(\phi_{k_1}^2 (x) \phi_{k_2}^2 (y)) = 1$. 
		\item $Corr(c_{k_1,k_2}, c_{l_1,l_2})=0$ if $(k_1,k_2) \neq (l_1,l_2)$.
	\end{enumerate}
\end{theorem}
The first condition in this theorem comes from the fact that all the basis functions are normalized. 
It also shows the variance of the frequency domain data for the copula-transformed independent data set. 
The second item shows how the coefficients
are correlated.  Note that there are many such necessary conditions for the frequency domain data. For example, 
when $N >> K$ and $X, Y \sim U[0,1]$ are independent, then each entry in the frequency domain coefficient 
matrix $C$ approximately follows a Gaussian distribution. Borrowing results from random matrix theory, the empirical 
distribution of the eigenvalues of matrix $C$ (Wigner matrix) will converge to certain limiting distributions. 
Interested readers are referred to \cite{anderson2010introduction,mehta2004random} and references 
therein for further results from random matrix theory.

Utilizing the computed frequency domain coefficient matrix $\hat{C}=\{ \hat{c}_{k_1,k_2} \}$ and results in 
Theorem~\ref{theorem:ind03}, we present two algorithms (statistics) for independence testing using the frequency
domain data. 

\vspace{0.1in}
{\noindent \bf Algorithm 1 (Max-FiT):} In the first algorithm, we focus on a statistic $z$ 
defined by $z=\|\hat{C}\|_{max}=\max_{k_1,k_2} |\hat{c}_{k_1,k_2}|$. When $X, Y \sim U[0,1]$ are independent
and $N$ is large, the new random variable (statistic) $z$ asymptotically follows the generalized extreme value (GEV)
distribution~\cite{de2006extreme,michailidis2007extreme,weisstein2021extreme}, a well-studied topic 
in statistics; its probability density function is available asymptotically and numerically. 
Algorithm 1 determines whether $X$ and $Y$ are independent by comparing $z$ with the 
corresponding independent case GEV distribution.

We refer to Algorithm 1 as the Maximal coefficient Frequency-domain independence Test (Max-FiT) 
and its analysis and properties are very similar to those of the maximal coefficient BET 
technique discussed in \cite{zhang2019bet}, where the Walsh basis is used. 

\vspace{0.1in}
{\noindent \bf Algorithm 2 (SVD-FiT):} In our second algorithm, we use the statistic $z=||\hat{C}||_2$, 
the largest singular value of the computed coefficient matrix $\hat{C}$,  or equivalently, we can study 
the largest eigenvalue of the Wishart (Laguerre) ensemble of a positive definite 
matrix $\hat{C}^T \cdot \hat{C}$. The random variable $z$ follows the Tracy-Widom 
distribution \cite{bai2010spectral,tao2012topics,tracy2002distribution} 
which is also a well-studied topic. The algorithm determines whether $X$ and $Y$ are 
independent by comparing $z$ with the Tracy-Widom distribution for independent data.
We refer to the second algorithm as the Singular Value Decomposition Frequency-domain 
independence Test (SVD-FiT). This test is closely related to the physical domain 
$\chi^2$-test.

In both algorithm implementations, as we have applied the copula transform to each variable, 
the GEV or Tracy-Widom distributions for the independent case can be pre-constructed to avoid the 
expensive permutation tests. One can also apply approximate distributions that are asymptotically 
correct for large sample sizes. Note that in two dimensions, a general data set  
$\{ (x_n,y_n) \}_{n=1}^N$ can be efficiently transformed into a copula data set by applying 
a quick sort to each of the two random variable sets $\{ x_n \}$ and $\{y_n \}$, $n=1, \cdots, N$,
respectively. However, for higher-dimensional data, e.g., when $x_n$ is not a scalar but a vector 
or tensor, ``optimally" transforming arbitrary higher-dimensional data into a high-dimensional uniform distribution is still an active research area (e.g., in 
optimal transport). Therefore, in more general higher-dimensional settings, permutation tests 
or asymptotic analyses have to be applied to estimate the distributions for the independent cases.

\subsection{Finding Low-Dimensional Structures}
\label{sec:bump}
It is possible and sometimes numerically more convenient (in terms of stability and efficiency) 
to extract the hidden features of physical domain data from the frequency domain coefficient
matrix and higher order moments. 
In this section, we first discuss how to find the ``bump" functions defined as 
$y=H(x)=\Eb [Y|X = x]$ or $x=G(y)=\Eb [X|Y = y]$. We then generalize the ideas to find the 
low-dimensional structures (one dimensional curves in two dimensions) given by the implicit 
algebraic relation $H(x)=G(y)$ where the data are the perturbed random variables
around the curve implicitly given by $H(x)=G(y)$, e.g., the unit circle $x^2 = 1 - y^2$. 
As we will demonstrate in \cref{sec:num},
finding the separable relation $H(x)=G(y)$ is already numerically challenging and designing
a numerically stable algorithm requires better knowledge of many factors, including the 
underlying noise. Developing numerically stable formulations for more general relations $H(x,y)=0$ 
is currently being studied, and results will be reported in a subsequent paper.
 
\subsubsection{Bump Hunting}
Bump hunting is a well-studied topic in statistics and applications \cite{friedman1999bump}. 
In bump hunting, for a given set of physical domain data $(x_n,y_n)$, one searches for a 
function $H(x)$ (or $G(y)$), such that for any fixed $x$ (or $y$), $\Eb [Y - H(X)|X = x] = 0$ 
(or $\Eb [X - G(Y)|Y = y] = 0$). Using the joint probability density function $f(x,y)$
and its expansion $f(x,y)= \sum_{i,j} c_{i,j} \phi_i(x) \phi_j(y)$, the bump function satisfies
\begin{equation}
\label{eq:bump1}
 \int (y - H(x)) \left(\sum_{i,j} c_{i,j} \phi_i(x) \phi_j(y)\right)dy = 0 \mbox{ for any $x$,} 
\end{equation}
or 
\begin{equation}
\label{eq:bump2}
  \int (x - G(y)) \left(\sum_{i,j} c_{i,j} \phi_i(x) \phi_j(y)\right)dx = 0 \mbox{ for any $y$.} 
\end{equation}
Using $H(x)$ in Eq.~(\ref{eq:bump1}) as an example, as
$$\int y\sum_{i,j} c_{i,j} \phi_i(x) \phi_j(y)dy = \int H(x) \sum_{i,j} c_{i,j} \phi_i(x) \phi_j(y)dy$$
and the basis functions are orthonormal with $\phi_0=\frac{1}{||1||}$, we have 
$$\sum_{i} \left( \sum_{j} c_{i,j} \langle y, \phi_j(y)\rangle \right) \phi_i(x) = H(x) \langle1, \phi_0(y)\rangle 
	 \sum_{i} c_{i,0}  \phi_i(x), $$
where the inner product $\langle f(x),g(x)\rangle = \int f(x) g(x) dx$. Therefore,
 $$ H(x) = \frac{\sum_{i} \left( \sum_{j} c_{i,j} \langle y, \phi_j(y)\rangle \right)  \phi_i(x)}{\langle 1, \phi_0(y)\rangle 
	\sum_{i} c_{i,0}  \phi_i(x)}.$$
When the Legendre polynomial basis is used and $y \in [-1,1]$, the formula can be further simplified to 
 $$ H(x) = \frac{\sum_{i} \left( c_{i,1} \langle y, \phi_1(y)\rangle \right)  \phi_i(x)}{\langle 1, \phi_0(y)\rangle 
	\sum_{i} c_{i,0}  \phi_i(x)}$$
as $y$ is one of the (scaled) basis functions.

For two-dimensional implementations, the data set  $\{ (x_n,y_n) \}_{n=1}^N$ can be 
efficiently transformed into a copula data set. The transformed data set satisfies the 
condition that the marginal distribution of $X$ follows the uniform distribution $U[0,1]$ and 
$\sum_{i} c_{i,0}  \phi_i(x)=1$. Numerically, the copula transformation is preferred for
better numerical stability (and accuracy) properties, as when using the estimated coefficients $\hat{c}_{i,j}$,
the function $\sum_{i} \hat{c}_{i,0}  \phi_i(x)$ may not always be positive, and small perturbations
of this function may cause large errors in the resulting $H(x)$ estimate.  
We present a stable bump-hunting formulation when using the Legendre polynomial basis for 
a copula physical domain data set.

\vspace{0.1in}
{\noindent \bf Algorithm 3 (Bump Hunting):} Assume the data set 
$\{ (x_n,y_n) \}_{n=1}^N$ is a copula data set in $[0,1]^2$. One can 
construct the random coefficient matrix $C$ using the (normalized) 
Legendre polynomial basis (efficiently using fast Legendre transforms) 
and the bump function $H(x)$ is given by 
$$ H(x) = \frac{\sum_{i} \left( c_{i,1} \langle y, \phi_1(y)\rangle \right)  \phi_i(x)}{\langle 1, \phi_0(y)\rangle},$$
For more general basis choices with $\phi_0(x)=\frac{1}{\|1\|}$, the bump function 
is given by 
$$H(x) = \frac{\sum_{i} \left( \sum_{j} c_{i,j} \langle y, \phi_j(y)\rangle \right)  \phi_i(x)}{\langle 1, \phi_0(y)\rangle}.$$
The formulas for $G(y)$ can be derived similarly and we skip the details.

\subsubsection{Low Dimensional Structures $H(x)=G(y)$}
When the physical domain data $(x,y)$ are generated using the equation 
$y=\sin(x) + \epsilon$ where $\epsilon \sim N(0, \delta)$, bump hunting can effectively reveal 
the deterministic relationship $y=\sin(x)$. However, when the data pairs are generated around 
a circle $x^2+y^2=1$ with some independent random perturbations to $X$ and $Y$, respectively, 
because of the nonlinearity in the low-dimensional geometric structure (circle), the bump 
functions $\Eb [Y|X = x]$ and $\Eb [X|Y = y]$ are approximately given by the perpendicular lines 
$x=0$ and $y=0$, respectively. Neither of these corresponds to the low-dimensional structure we are searching for. 

Finding general low-dimensional structures hidden in the multi-variable probability density
function is an open problem and requires assumptions on many factors, including the underlying
geometric structure and properties of the random noise. 
For the one-dimensional structure $H(x)=G(y)$ in a two-dimensional data set, an intuitive 
idea is to search for the functions $H(x)$ and $G(y)$ such that 
$$\Eb [G(Y) - H(X)|X = x] = 0 \quad \mathrm{and} \quad \Eb [G(Y) - H(X)|Y = y] = 0.$$ 
Or more generally, $\Eb [H(X)-G(Y)|(\cos(\theta)X+\sin(\theta)Y =  
\cos(\theta) x +\sin(\theta) y)] =0$ for all possible angles $\theta$. This is possible
when data are exactly located on the low-dimensional manifold and there is no noise. 
However, for noisy data, there are existence, uniqueness, regularity issues, as well as 
computational efficiency and stability issues with this formulation. For example, consider 
$x=1+\epsilon$ where $\epsilon \sim N(0, \delta)$ and the deterministic relation is given by $x=1$. If one studies $x^2$, a simple calculation shows $\Eb (x^2-1) \neq 0$ and
one gets a (slightly) different relation. Another example is $y=\sin(x) + \epsilon$, where 
$-1 \leq x \leq 1$ and $\epsilon \sim N(0, \delta)$. The deterministic relation using
$\Eb [Y|X=x]$ gives the function $y=\sin(x)$, while bump hunting returns a (slightly) different
relation if we use $\Eb (X|Y = y)$. 

In this paper, we introduce a specific definition of the low-dimensional structure given by
the relation $H(x)=G(y)$, which is well-defined and efficient to calculate. We first present
the algorithms (based on existing statistical data analysis tools).

{\noindent \bf Algorithm 4.1 (CCA - Low Dimensional Structure):} For a given data set 
$\{ (x_n,y_n) \}_{n=1}^N$, (not necessarily a copula data set in $[0,1]^2$),  
construct the frequency domain data sets 
$$X=\{ \phi_0(x_i), \phi_1(x_i), \cdots, \phi_K(x_i) \} \quad \mathrm{and} \quad  
Y=\{ \phi_0(y_i), \phi_1(y_i), \cdots, \phi_K(y_i) \}, \quad i=1, \cdots, N.$$ 
Apply existing canonical correlation analysis (CCA) tools, for example, the Matlab function 
$$[\vec A,\vec B,\vec{r},\vec{U},\vec{V}] = canoncorr(X,Y).$$ Then the CCA-low-dimensional structure (CCA-LDS) is given by $H(x)=G(y)$ where
$$H(x)= \sum_{k=0}^K A(k,1) \phi_k(x), \quad \quad G(y)= \lambda \sum_{k=0}^K B(k,1) \phi_k(y),$$ 
and $\lambda$ is a constant which can be computed using a simple linear fit.

Note that most existing CCA tools do not consider the special structures of the basis
functions. As the CCA algorithm is based on the singular value decomposition algorithm,
to take advantage of existing fast transformation algorithms (e.g., NUFFT  
or fast Haar transform), we present a slightly modified algorithm for finding
$A(:,1)$ and $B(:,1)$.

{\noindent \bf Algorithm 4.2 (SVD - Low Dimensional Structure):} For a given data set 
$\{ (x_n,y_n) \}_{n=1}^N$, (not necessarily a copula data set in $[0,1]^2$),  
construct the frequency domain coefficient and variance matrices $\hat{C}$, $\tilde{\Sigma}_{XX}$, 
$\tilde{\Sigma}_{YY}$, and $\tilde{\Sigma}_{XY}$ using a fast transformation algorithm.
Apply existing singular value decomposition (SVD) tools to the covariance matrix $\tilde{\Sigma}_{XY}$,
and denote the left and right singular vectors corresponding to the largest singular value
as $\vec{c}$ and $\vec{d}$. Define $\vec{\omega}=\tilde{\Sigma}_{XX}^{-1/2} \vec{c}$ and 
$\vec{\beta}=\tilde{\Sigma}_{YY}^{-1/2} \vec{d}$.
Then the SVD-low-dimensional structure (SVD-LDS) is given by $H(x)=G(y)$ where
\begin{equation}
	\label{eq:hgexpan}
	H(x)= \sum_{k=1}^K \omega_k \phi_k(x), \quad \quad G(y)= \lambda \sum_{k=1}^K \beta_k \phi_k(y) + \gamma,
\end{equation}
and $\lambda$ and $\gamma$ are constants that can be computed using a simple linear fit.

{\noindent \bf Comment:} Note that the index $k$ starts from $0$ in Algorithm 4.1, but from $1$ in Algorithm 4.2.
Also, Algorithm 4.2 can be more efficient due to the application of fast transformation algorithms. 
In Algorithms 4.1 and 4.2, the coefficients $A(k,1)$, $B(k,1)$, 
$\omega_k$, and $\beta_k$ are uniquely determined. However, the coefficients $\lambda$ 
and $\gamma$ depend on how the linear fit problem is formulated. In the current implementation,
we solve the least squares optimization problem
$$ \sum_{n=1}^N \left( \sum_{k=1}^K \omega_k \phi_k(x_n) - 
           (\lambda \sum_{k=1}^K \beta_k \phi_k(y_n) + \gamma ) \right)^2.$$
An alternative least squares problem is 
$$  \sum_{n=1}^N \left( \tilde{\lambda} \sum_{k=1}^K \omega_k \phi_k(x_n) +\tilde{\gamma} 
                             - \sum_{k=1}^K \beta_k \phi_k(y_n) \right)^2$$
and the functions $H(x)$ and $G(y)$ are then defined as
$$	H(x)= \tilde{\lambda} \sum_{k=1}^K \omega_k \phi_k(x_n) +\tilde{\gamma} 
            \quad \mathrm{and} \quad G(y)= \sum_{k=1}^K \beta_k \phi_k(y).$$
As we will demonstrate in \cref{sec:num}, the two different formulations 
define different low-dimensional manifold structures. Indeed, the proper definition 
of physically relevant low-dimensional structures in the data requires knowledge
of the structures and assumptions of the data set and noise.

\subsubsection{Algebraic Interpretation} We explain the algebraic/geometric meaning of the relation equation $H(x)=G(y)$ 
and its connections to the bump hunting type conditions 
\begin{equation}
	\label{eq:2conditions}
	\Eb [H(X)-G(Y)|X = x]=0 \mbox{ and } \Eb [H(X)-G(Y)|Y = y]=0.
\end{equation} 
As there are no $\{ H(x), G(y) \}$ that exactly satisfy both conditions in 
Eq.~(\ref{eq:2conditions}), in the remainder of this section, we show that 
the relation $H(x)=G(y)$ from Algorithm 4.2 (the same as that from Algorithm 4.1) 
represents a particular ``best" least squares type solution to the inconsistent 
system in Eq.~(\ref{eq:2conditions}).

We first consider a fixed $x$ and the condition $\Eb [G(Y) - H(X)|X = x] = 0$. We start with
$$ \int (G(y) - H(x))(\sum_{i,j} c_{i,j} \phi_{i} (x) \phi_{j} (y)) dy=0,$$
or equivalently,
$$\int G(y)(\sum_{i,j} c_{i,j} \phi_{i} (x) \phi_{j} (y)) dy = \int H(x)(\sum_{i,j} c_{i,j} \phi_{i} (x) \phi_{j} (y)) dy = H(x) P_X(x), $$
where $P_X(x)$ is the marginal probability density function of the random variable $X$.
When $X$ $\sim U[0,1]$, $P_X(x)=1$.
Next, using the expansions 
$$ H(x)= \sum_{k=0}^K \omega_k \phi_k(x) \quad \mathrm{and} \quad G(y)= \sum_{k=0}^K \beta_k \phi_k(y)$$
and the orthogonality of the 
basis functions, we have the equations for the coefficients
$$ C  \left[ \begin{array}{c} \beta_0 \\ \beta_1 \\ \vdots \\ \beta_K \end{array}\right]  = \left(
	\Sigma_{XX} + \left[ \begin{array}{c} c_{0,0} \\ c_{1,0} \\ \vdots \\ c_{K,0} \end{array} \right] 
	[ c_{0,0}, c_{1,0}, \cdots, c_{K,0} ] \right)  
	\left[ \begin{array}{c} \omega_0 \\ \omega_1 \\ \vdots \\ \omega_K \end{array}\right],
$$
where $C$ is the frequency domain coefficient matrix and $\Sigma_{XX}$ is the variance matrix
defined in \cref{sec:transformed}.
As $\phi_0(x)$ is the constant function, we further simplify the relation equation by 
introducing $\tilde{C}$
$$ C  = \left[ \begin{array}{cccc} c_{0,0} & c_{0,1} & \cdots & c_{0,K} \\
	c_{1,0} & & & \\
		\vdots & & \tilde{C} & \\
		c_{K,0} & & & \\
	\end{array} \right]
$$
and consider the new relation equation
$$\left[ \begin{array}{c} 
	\sum_{k=0}^K c_{0,k} \beta_k \\
	\beta_0  \left[ \begin{array}{c} c_{1,0} \\ \vdots \\ c_{K,0} \end{array} \right] + 
		\tilde{C}  \left[ \begin{array}{c} \beta_1 \\ \vdots \\ \beta_K \end{array} \right] \\
\end{array} \right]  = \left[ \begin{array}{c} 
	\sum_{k=0}^K c_{k,0} \omega_k \\
	\tilde{\Sigma}_{XX}  \left[ \begin{array}{c} \omega_1 \\ \vdots \\ \omega_K \end{array} \right] +
		(\sum_k c_{k,0} \omega_k)  \left[ \begin{array}{c} c_{1,0} \\ \vdots \\ c_{K,0} \end{array} \right] \\
\end{array} \right]. 
$$
Using the first equation $\sum_{k=0}^K c_{0,k} \beta_k = \sum_{k=0}^K c_{k,0} \omega_k$, we can
translate the second vector equation system into 
$$ (\beta_0 - \sum_{k=0}^K c_{0,k} \beta_k) \left[ \begin{array}{c} c_{1,0} \\ \vdots \\ c_{K,0} \end{array} \right] + 
		\tilde{C}  \left[ \begin{array}{c} \beta_1 \\ \vdots \\ \beta_K \end{array} \right] 
 = \tilde{\Sigma}_{XX}  \left[ \begin{array}{c} \omega_1 \\ \vdots \\ \omega_K \end{array} \right].
$$
Further simplifying the equation using the assumption that $\phi_0(x)$ is a constant function, we have
$$ (\tilde{C} - \left[ \begin{array}{c} c_{1,0} \\ \vdots \\ c_{K,0} \end{array} \right] [c_{0,1}, \cdots, c_{0,K}])  
	\left[ \begin{array}{c} \beta_1 \\ \vdots \\ \beta_K \end{array} \right] 
 = \tilde{\Sigma}_{XX}  \left[ \begin{array}{c} \omega_1 \\ \vdots \\ \omega_K \end{array} \right],
$$
or equivalently, 
\begin{equation}
\label{eq:lds01}
\tilde{\Sigma}_{XY} \left[ \begin{array}{c} \beta_1 \\ \vdots \\ \beta_K \end{array} \right] 
 = \tilde{\Sigma}_{XX}  \left[ \begin{array}{c} \omega_1 \\ \vdots \\ \omega_K \end{array} \right].
\end{equation}

Similarly, if we fix $y$ and consider the equation $\Eb [\tilde{G}(Y) - \tilde{H}(X)|Y = y] = 0$,
we get the linear system
\begin{equation}
\label{eq:lds02}
 \tilde{\Sigma}_{XY}^T \left[ \begin{array}{c} \tilde{\omega}_1 \\ \vdots \\ \tilde{\omega}_K \end{array} \right] 
 = \tilde{\Sigma}_{YY}  \left[ \begin{array}{c} \tilde{\beta}_1 \\ \vdots \\ \tilde{\beta}_K \end{array} \right].
\end{equation}
We used the notation $\tilde{H}$, $\tilde{G}$, $\tilde{\omega}$ and $\tilde{\beta}$ to represent a new set of 
functions and vectors. If we require the new set to be the same as the original ones, we will have 
$$ \left( \tilde{\Sigma}_{YY}^{-1}  \tilde{\Sigma}_{XY}^T \tilde{\Sigma}_{XX}^{-1}  \tilde{\Sigma}_{XY} \right)
	\left[ \begin{array}{c} \tilde{\beta}_1 \\ \vdots \\ \tilde{\beta}_K \end{array} \right]
 = \left[ \begin{array}{c} \tilde{\beta}_1 \\ \vdots \\ \tilde{\beta}_K \end{array} \right]
$$
and
$$ \left( \tilde{\Sigma}_{XX}^{-1}  \tilde{\Sigma}_{XY} \tilde{\Sigma}_{YY}^{-1}  \tilde{\Sigma}_{XY}^T \right)
	\left[ \begin{array}{c} \tilde{\omega}_1 \\ \vdots \\ \tilde{\omega}_K \end{array} \right] 
 = \left[ \begin{array}{c} \tilde{\omega}_1 \\ \vdots \\ \tilde{\omega}_K \end{array} \right],
$$
i.e., $1$ is an eigenvalue of the matrix 
$\tilde{\Sigma}_{XX}^{-1}  \tilde{\Sigma}_{XY} \tilde{\Sigma}_{YY}^{-1}  \tilde{\Sigma}_{XY}^T$. Such 
solutions unfortunately may not exist.
Instead of finding $\tilde{H}=H$, $\tilde{G}=G$, $\tilde{\omega}=\omega$, and $\tilde{\beta}=\beta$,
SVD-LDS and CCA-LDS solve an optimization problem. In the frequency domain, consider a 
given set of coefficients $[\beta_1, \beta_2, \cdots, \beta_K]$ (a given function $G(y)$). Define 
$$\vec{c} = {\Sigma_{YY}}^{\frac{1}{2}}\begin{bmatrix} {\beta}_1\\ \vdots\\ {\beta}_K\end{bmatrix}$$
and assume $||\vec{c}||=1$ (normalized $G(y)$ using a particular norm). One can compute the coefficient
vector $\vec{d} =  {\Sigma_{XX}}^{\frac{1}{2}}\begin{bmatrix} \omega_1\\ \vdots\\ \omega_K\end{bmatrix}$
(a resulting $H(x)$ function) using Eq.~(\ref{eq:lds01}) to get 
$$\vec{d}= \Sigma_{XX}^{-1/2}\Sigma_{XY} \Sigma_{YY}^{-\frac{1}{2}}\vec{c}.$$ 
Then using $\vec{d}$ and Eq.~(\ref{eq:lds02}), one can get a new coefficient vector
$$\tilde{c} = \Sigma_{YY}^{-\frac{1}{2}}\Sigma_{YX}\Sigma_{XX}^{-\frac12} \vec{d}$$
(a new function $\tilde{H}$). We have the following theorem
\begin{theorem}
The algorithms 4.1 and 4.2 return a vector $\| \vec{c} \|=1$ which minimizes 
$$||\tilde{c}-\vec{c}|| = \|(\Sigma_{YY}^{-\frac{1}{2}}\Sigma_{YX}\Sigma_{XX}^{-1}\Sigma_{XY}\Sigma_{YY}^{-\frac{1}{2}} - I) \vec{c} \|.$$ 
\end{theorem}
A similar result can be derived when one starts from $\vec{d}$ and
minimizes $\|\tilde{d}-\vec{d}\|$.  Comparing the results with standard CCA analysis, we summarize
our results in the following theorem for the vectors $\vec{\beta}$, $\vec{\omega}$, $\vec{c}$, and $\vec{d}$.
\begin{theorem}
The solution vectors $\vec{\beta}$, $\vec{\omega}$, $\vec{c}$, and $\vec{d}$ from Algorithms
4.1 and 4.2 have the following properties. 

{\noindent 1.} The vectors $\vec{\beta}$ and $\vec{\omega}$ maximize  
	$$\rho(\vec{\beta},\vec{\omega})= \frac{\vec{\beta}^T \tilde{\Sigma}_{XY} \vec{\omega}}
  {\sqrt{\vec{\beta}^T \tilde{\Sigma}_{XX} \vec{\beta}} \sqrt{\vec{\omega}^T \tilde{\Sigma}_{YY} \vec{\omega}} },$$
or in statistical language,
$$\operatornamewithlimits{argmax}\limits_{\vec{\beta},\vec{\omega}} Corr( \vec{\beta}^T \phi(x), \vec{\omega}^T\phi(y))
 = \operatornamewithlimits{argmax}\limits_{\vec{\beta},\vec{\omega}} \frac{E[(\vec{\beta}^T \phi(x))(\vec{\omega}^T\phi(y))]}
{\sqrt{E[( \vec{\beta}^T \phi(x))^{2}]E[(\vec{\omega}^T\phi(y))^{2}]}}. $$

{\noindent 2. } The vectors $\vec{c}$ and $\vec{d}$ maximize  
	$$\rho= \frac{\vec{c}^T  \Sigma_{XX}^{-\frac{1}{2}}\tilde{\Sigma}_{XY} \Sigma_{YY}^{-\frac12} \vec{d} }
  {\sqrt{\vec{c}^T \vec{c}} \sqrt{\vec{d}^T \vec{d}} },$$
	where $\vec{c} = {\Sigma_{YY}}^{\frac{1}{2}} \vec{\beta}$ and  $\vec{d}= {\Sigma_{XX}}^{\frac{1}{2}} \vec{\omega}$.

{\noindent 3. } $\vec{c}$ and $\vec{d}$ are the singular vectors of the correlation matrix $\Sigma_{XY}$. 

{\noindent 4. } $\vec{c}$ and $\Sigma_{XX}^{-\frac{1}{2}}\tilde{\Sigma}_{XY} \Sigma_{YY}^{-\frac12} \vec{d}$ are collinear,
	$\vec{d}$ and $\Sigma_{YY}^{-\frac{1}{2}}\tilde{\Sigma}_{YX} \Sigma_{XX}^{-\frac12} \vec{c}$ are collinear. 
\end{theorem}
As we will demonstrate in \cref{sec:num}, Algorithms 4.1 and 4.2 can be applied to find specific 
low-dimensional structures described by $H(x)=G(y)$, especially when both $H(x)$ and $G(y)$ only contain low-frequency modes. However, we found that the current formulations are not the most stable numerically. A more serious problem is that although the new nonlinear 
relation $H(x)=G(y)$ is more general compared with the relations $y=f(x)$ and $x=g(y)$ in 
the bump hunting algorithms, some relations cannot be described by $H(x)=G(y)$. For
example, none of the proposed algorithms can effectively find the (polynomial type) 
relation equation $x \cdot y =0$ or more general algebraic varieties. For more general 
low-dimensional structures described by the equation $H(x,y)=0$, the theoretical analysis 
and efficient and stable numerical solutions of a well-defined optimization problem are 
currently being studied, and results will be presented in the future.

%% file: sec04Numerical_0503_Juan.tex
\section{Preliminary Numerical Results}
\label{sec:num}
In this section, we present some preliminary numerical results for the frequency domain-based statistical 
data analysis algorithms.

\subsection{Dependency Test}
\begin{figure}[htbp]
    \centering
    \includegraphics[height=7cm,width=15.5cm]{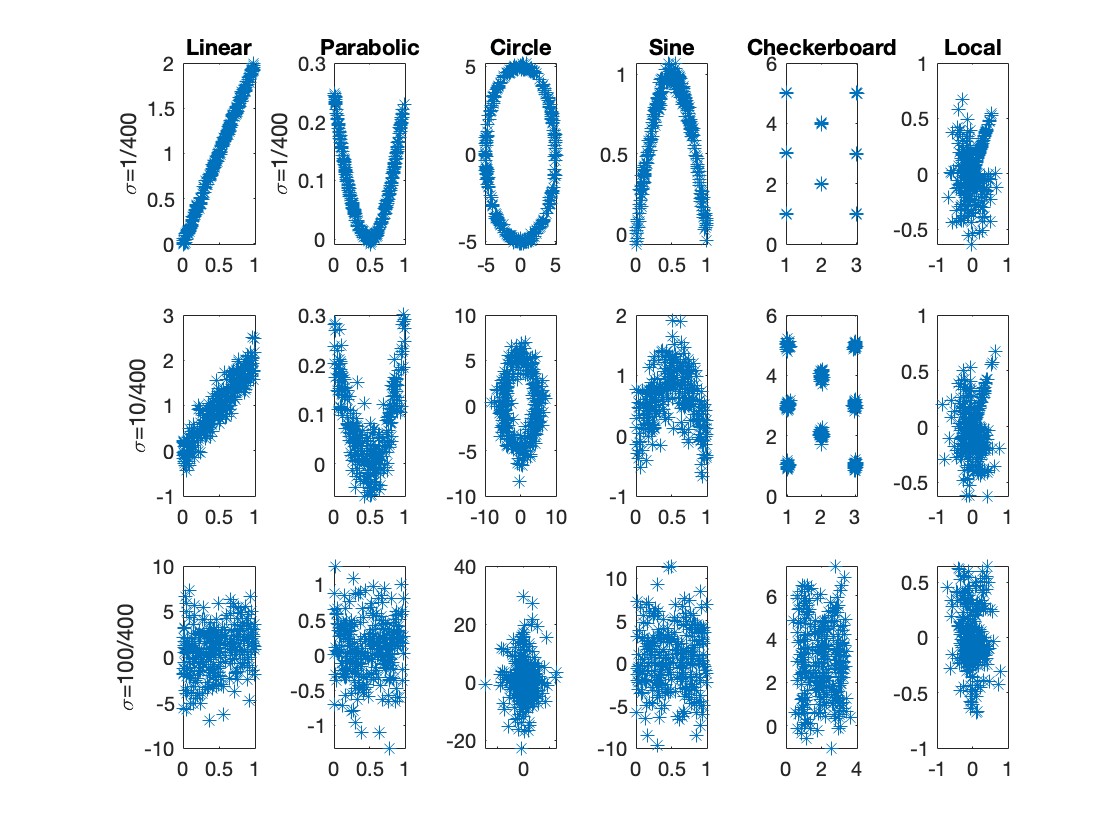}
	\caption{Six dependent random data sets for different $\sigma$ values $1/400$ (row 1), 
	$10/400$ (row 2), and $100/400$ (row 3).
    \label{fig:examples} }
\end{figure}

We demonstrate the statistical power of the frequency domain dependency test algorithms Max-FiT 
and SVD-FiT. We present the results for the Fourier series, Legendre polynomials and Haar wavelet basis 
functions, and compare them with those from the Walsh functions-based BET algorithms in~\cite{zhang2019bet} 
for 6 types of manufactured dependent data sets in Table~\ref{tab:examples}, 
where $U\sim Uniform[0, 1]$, $W \sim MultiBern (\{1, 2, 3\},(1/3, 1/3, 1/3))$,
$V_1\sim Bern(\{2, 4\},(1/2, 1/2))$, $V_2\sim MultiBern(\{1, 3, 5\},(1/3, 1/3, 1/3))$, 
$G_1, G_2 \sim N (0, 1/4)$, and
$\epsilon$, $\epsilon'$, and $\epsilon''$ $ \sim N(0,\sigma^2)$. 
%\lipsum[50]
\begin{table}[htbp]
	\footnotesize
	\caption{Manufactured dependent data sets for independence tests.} 
	\label{tab:examples}
\begin{center}
\begin{tabular}{|c|c|c|}
\hline
Linear & $X=U$ & $Y = 2X+10\epsilon$\\
\hline
%Name & X & Y \\
%\hline
%\hline
Parabolic & $X=U$ & $Y=(X-0.5)^2+1.5 \epsilon$ \\
\hline
Circular & $X=5\cos \theta + 40 \epsilon$ & $Y=5\sin \theta + 40 \epsilon'$ \\
\hline
Sine & $X=U$ & $Y=\sin(2 \pi X) +15 \epsilon$ \\
\hline
Checkerboard & $X=W+\epsilon$ & $Y= \left\{ \begin{array}{ll} V_1+ \epsilon' \quad & if \quad W=2 \\ V_2+\epsilon'' \quad & otherwise. \end{array} \right. $ \\
\hline
Local & $X=G_1$ & $Y= \left\{ \begin{array}{ll} X + \epsilon' \quad & if \quad 0 \leq G_1, G_2 \leq 1 \\ G_2 \quad & otherwise. \end{array} \right. $ \\
\hline
\end{tabular}
\end{center}
\end{table}
Sampled data sets for different $\sigma$ choices are plotted in Figure~\ref{fig:examples}. 

\begin{figure}[htbp]
    \centering
    \subfloat[$Linear$]{\label{fig:a}\includegraphics[width = 3in]{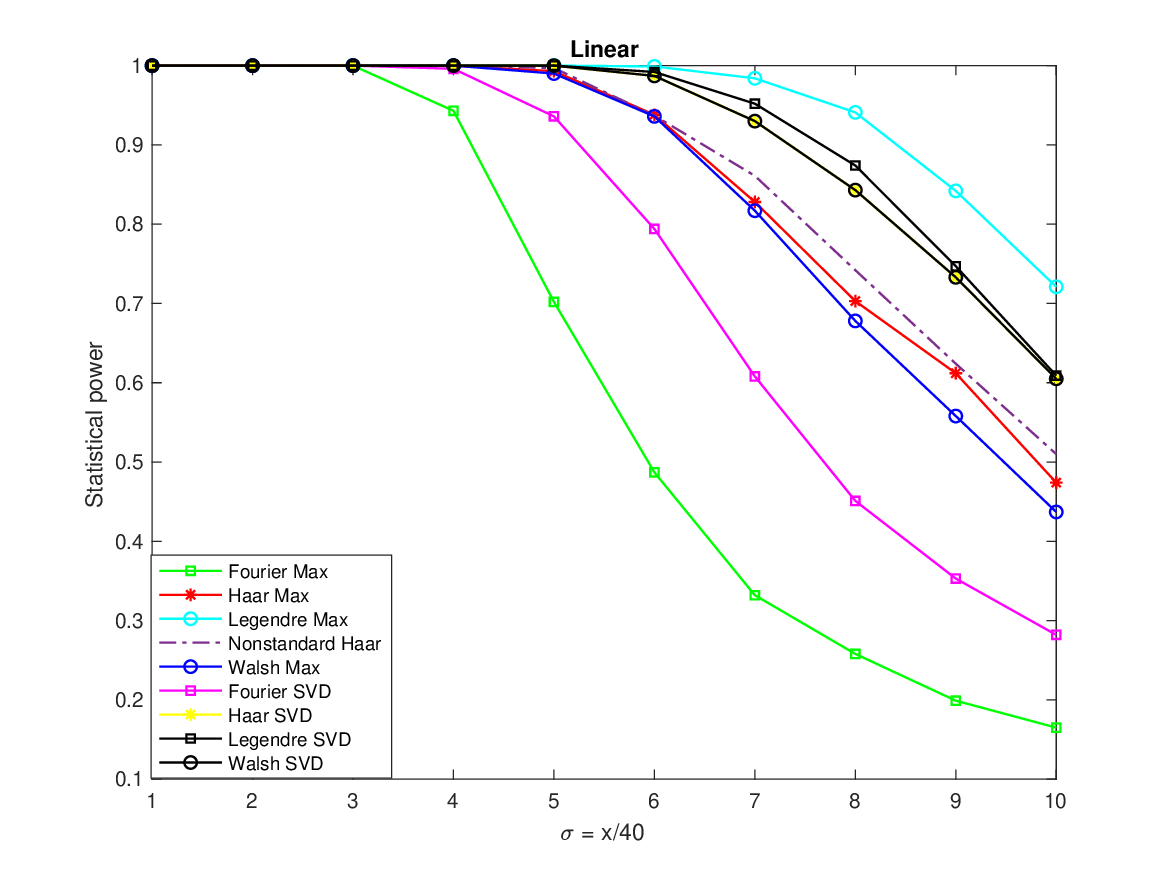}}
    \subfloat[$Parabolic$]{\label{fig:b}\includegraphics[width = 3in]{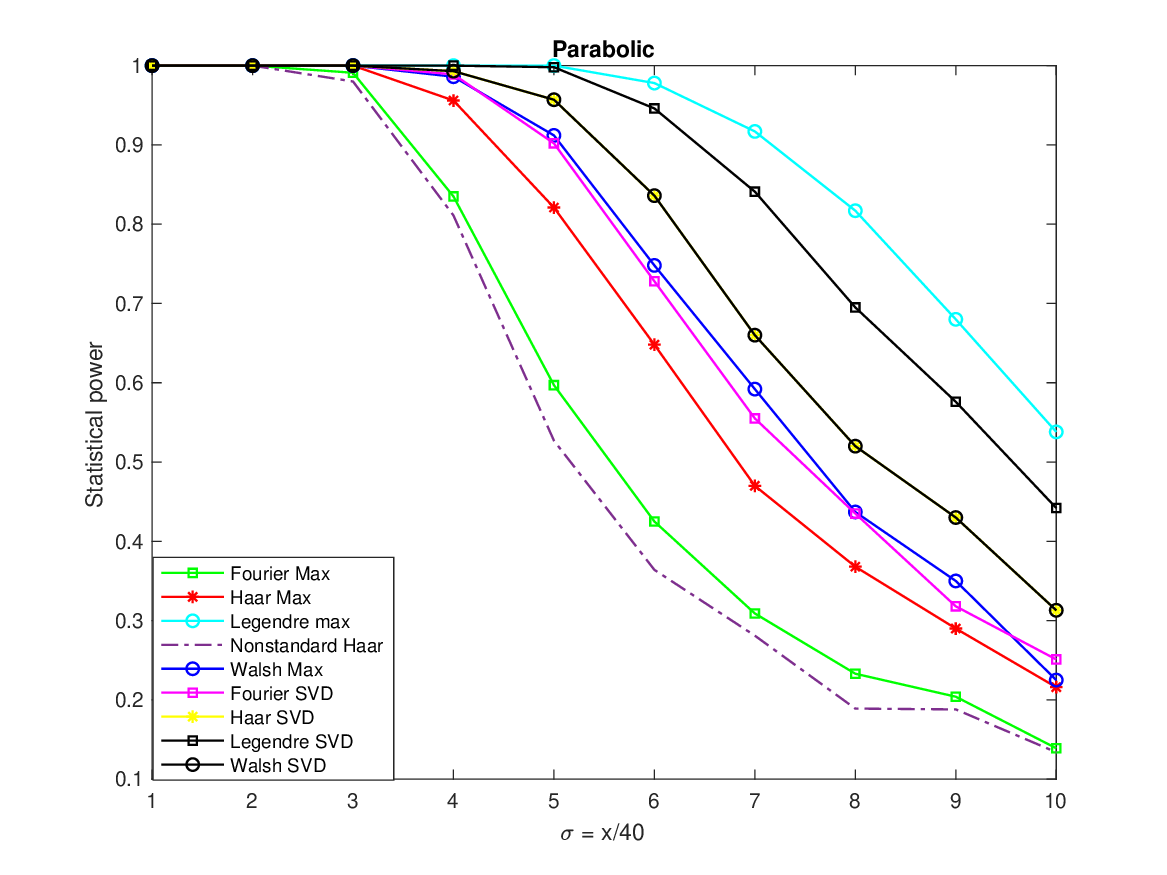}}\\
    \subfloat[$Circle$]{\label{fig:c}\includegraphics[width = 3in]{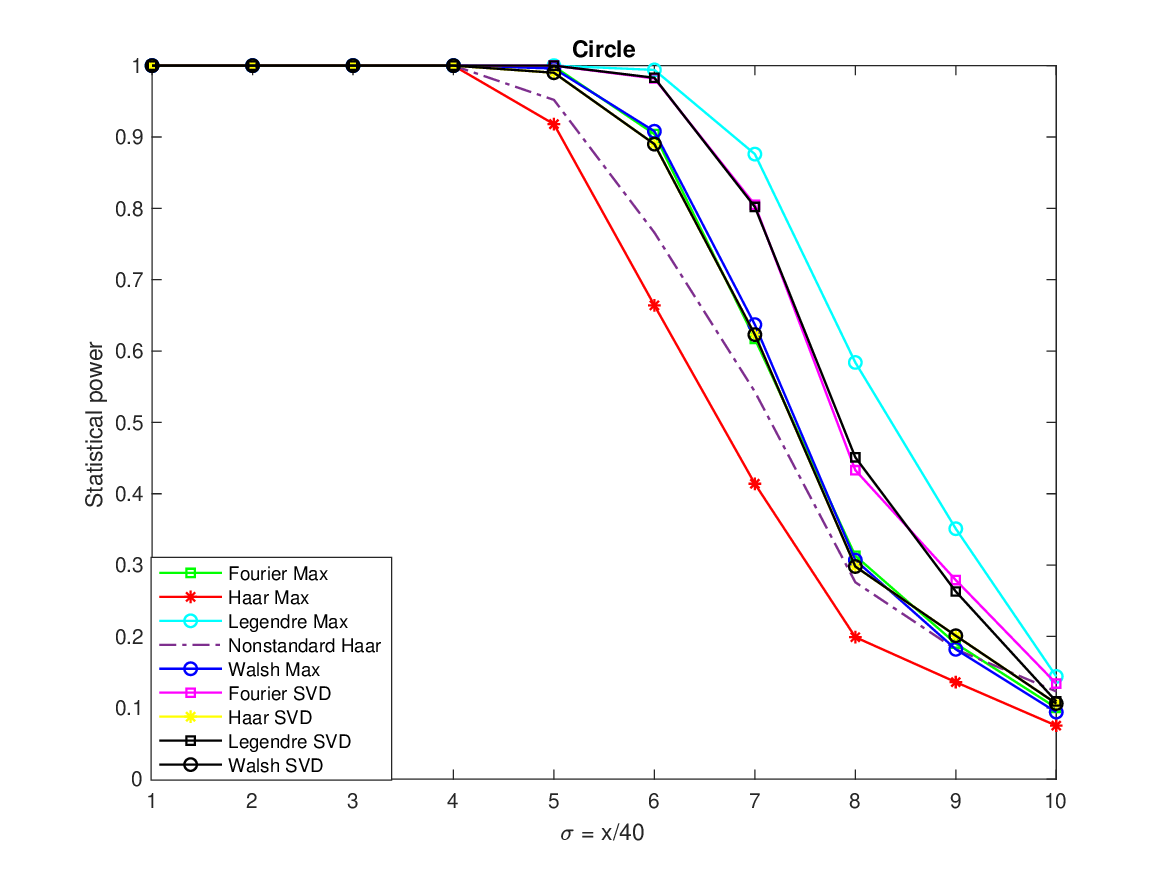}}
    \subfloat[$Sine$]{\label{fig:d}\includegraphics[width = 3in]{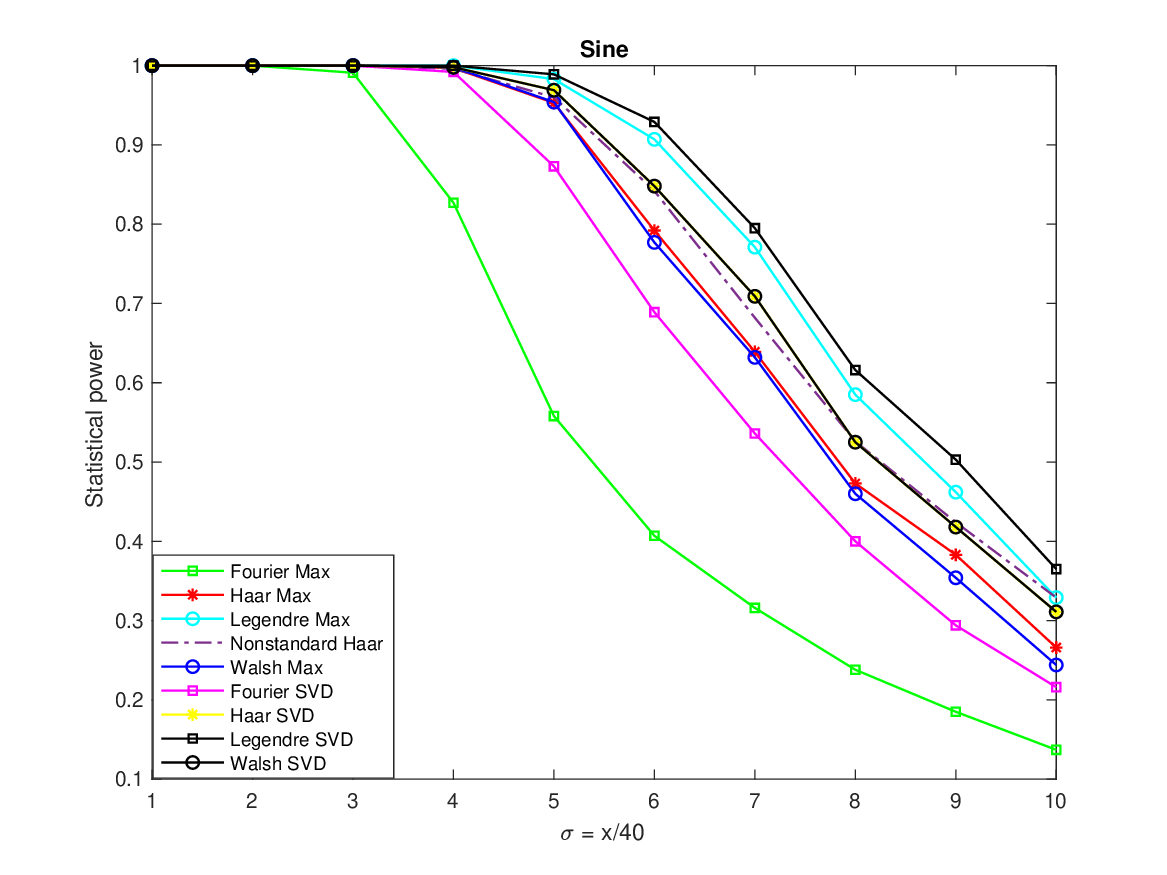}}\\
    \subfloat[$Checkerboard$]{\label{fig:e}\includegraphics[width = 3in]{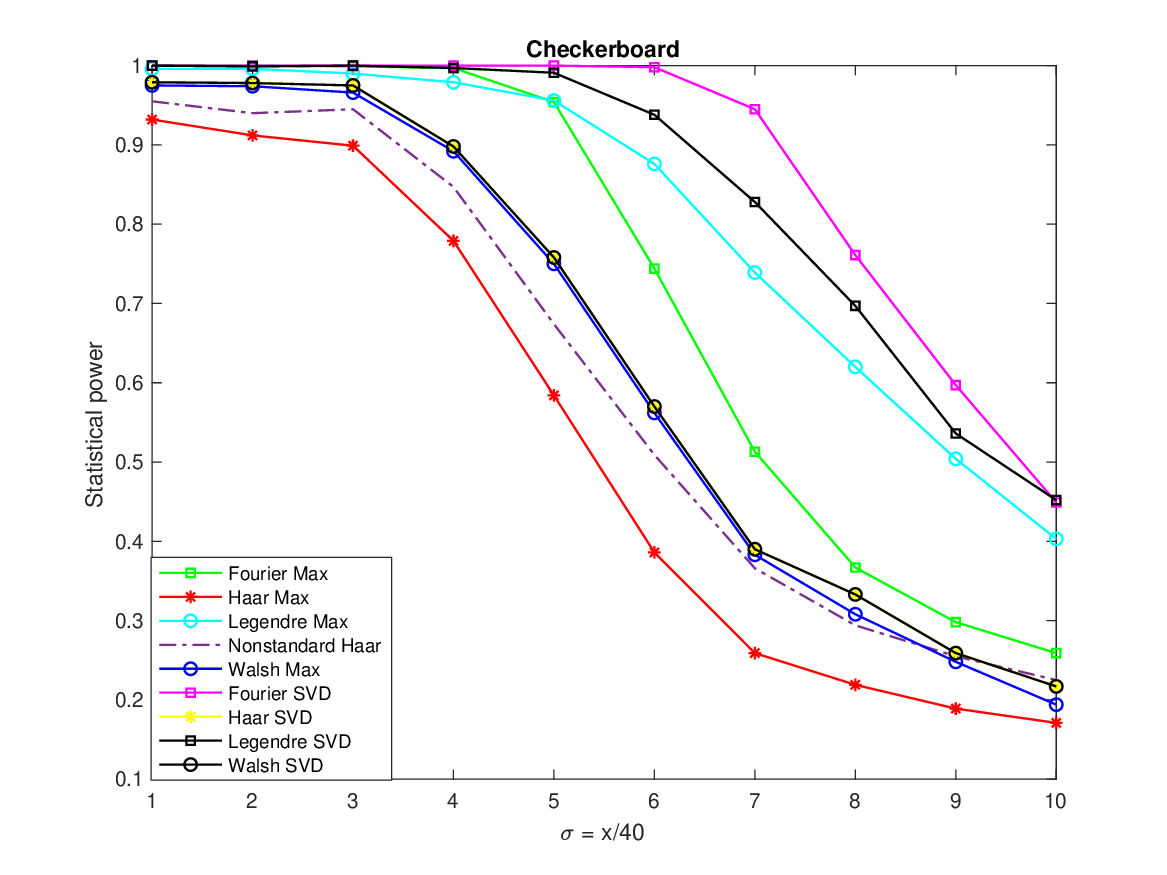}}
    \subfloat[$Local$]{\label{fig:f}\includegraphics[width = 3in]{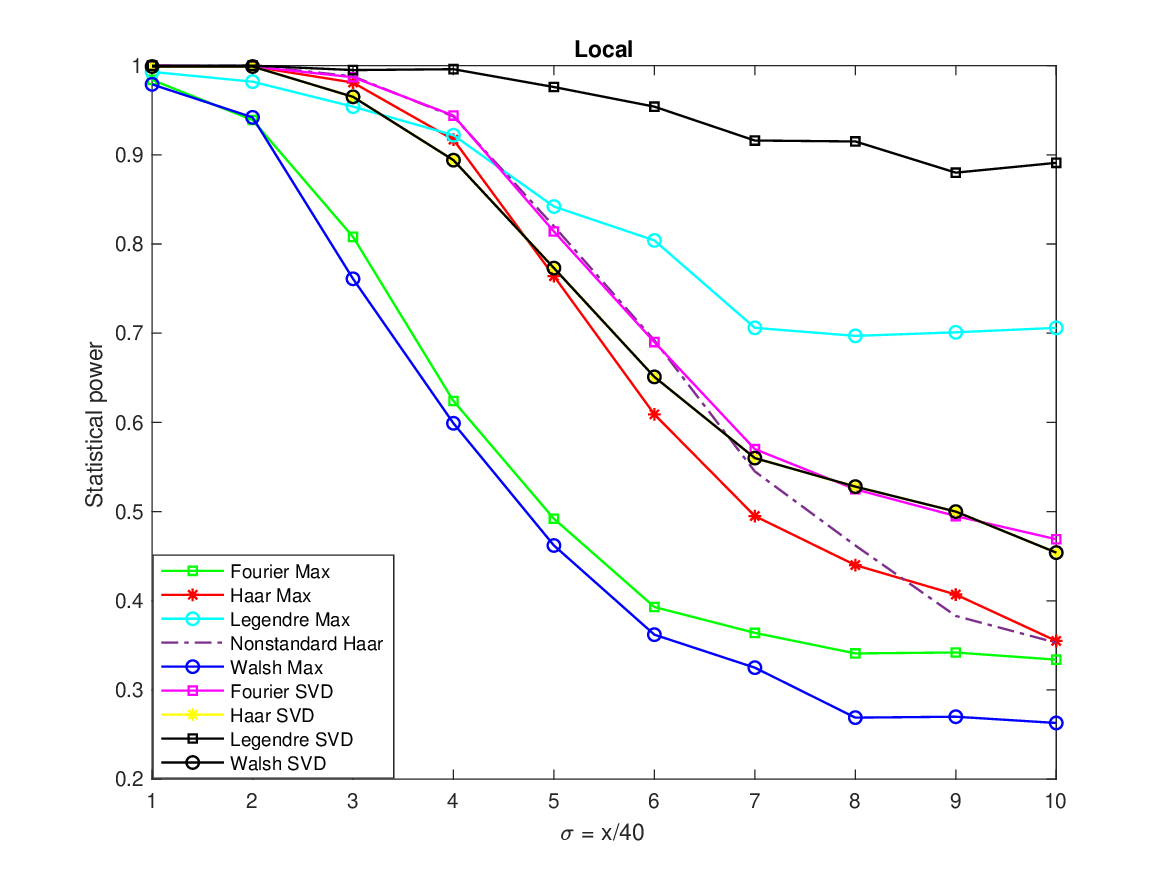}}\\
    \caption{Statistical power of Max-FiT and SVD-FiT for different bases.
    \label{fig:l_p}}
\end{figure}
In Figure~\ref{fig:l_p}, we show the statistical power of the Max-FiT and SVD-FiT algorithms
using different basis functions. The statistical power is the probability that the test
correctly rejects the independence hypothesis for a dependent data set. In the plot, 
the $x$-axis shows the noise level $\sigma=x/40$ and the $y$-axis shows the statistical
power. For each noise level, we generate $1000$ random data sets to estimate the 
probability. Each data set has $256$ data points $\{ (x_n,y_n) \}_{n=1}^{256}$. 
We simulate the null distributions of the independent cases for different basis choices
and statistics using a permutation test~\cite{annis2005permutation,collingridge2013primer,
edgington2007randomization,fisher1936design} with $200$ 
permutations for each case. Note that the computationally expensive permutation tests can be avoided 
in future implementations by analytically approximating the null distributions using 
asymptotic analysis.
%One current project is to construct these distributions either analytically or using asymptotic analysis to avoid the very expensive permutation test. 
For the manufactured examples, we found that the Legendre 
polynomial-based algorithms (``Legendre Max" and ``Legendre SVD") in general outperform 
methods using other bases.  Methods based on the local Haar basis functions often 
perform poorly (``Haar Max" and ``Nonstandard Haar"), and the statistical power of the 
global ``Haar SVD" method is identical to that of the ``Walsh SVD" algorithm, 
which suggests that an optimal independence test statistic needs to consider the
global properties of the manufactured data sets. 

\begin{figure}[htbp]
    \centering
    \includegraphics[height=5cm,width=8cm]{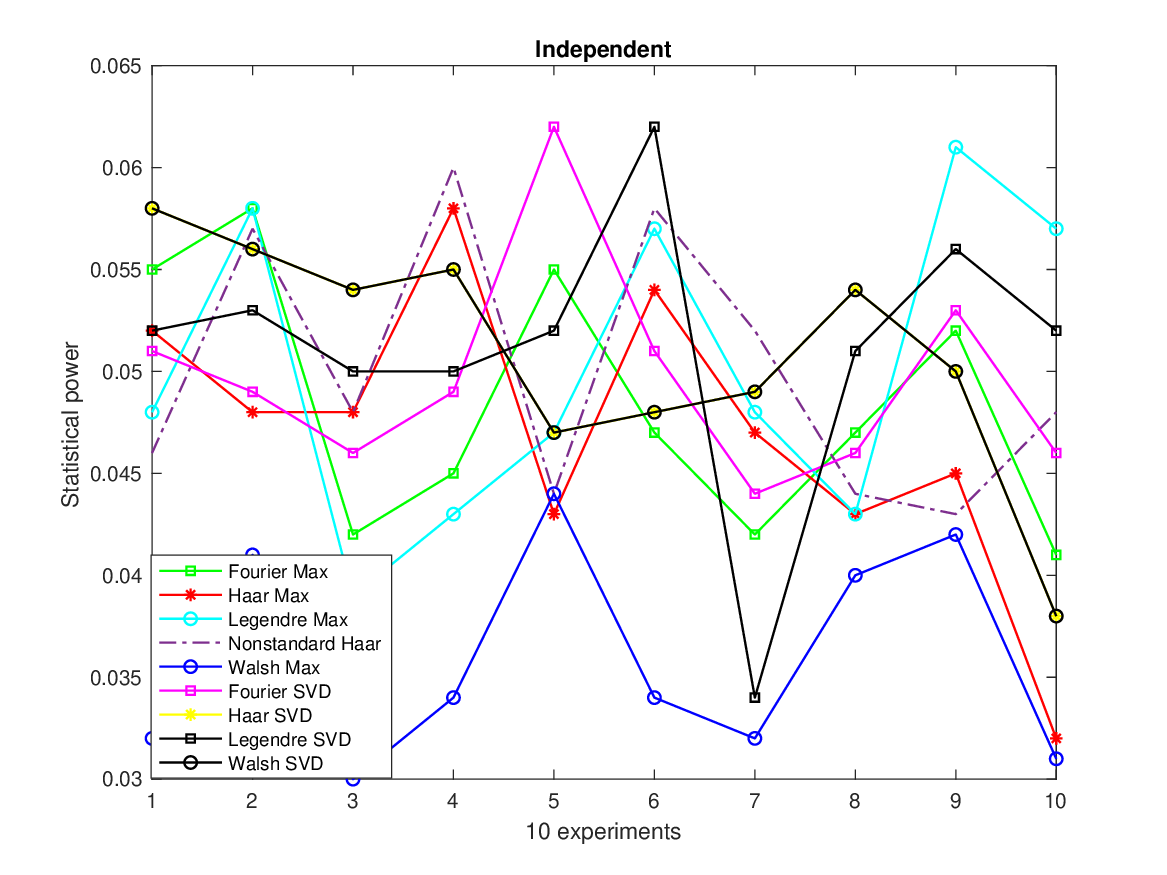}
	\caption{Statistical power of Max-FiT and SVD-FiT for independent data sets.}
    \label{fig:inddata}
\end{figure}

In Figure~\ref{fig:inddata}, we present the statistical power of the algorithms for independent
data sets. The experiments were repeated $10$ times ($x$-axis) and each experiment considered 
$1000$ random data sets each containing $256$ data points $\{ (x_n,y_n) \}_{n=1}^{256}$. 
The results assess the type I error of rejecting the independence hypothesis when
it is actually true. As we set the $p$-value threshold at $5\%$, the statistical
power of each algorithm ($y$-axis) correctly oscillates around $5\%$.

In addition to manufactured data sets, we also applied the algorithms to real 
world data sets. In Table~\ref{tab:stars}, we study the (dependent) galactic coordinates 
of the 256 brightest stars in the night sky from~\cite{zhang2019bet}. The galactic 
coordinates are essentially spherical coordinates with the Sun at the center.
\begin{table}[htbp]
\caption{p-values from different algorithms for galactic coordinates of 256 brightest stars.}
	\label{tab:stars}
\begin{center}
\begin{tabular}{|c|c|c|c|}
\hline
 & p-value &  & p-value\\
\hline
Fourier Max & 0.03 & Fourier SVD & 0.000\\
\hline
Haar Max & 0.025 & Haar SVD &0.005\\
\hline
 Legendre Max& 0.040  & Legendre SVD  &0.000\\
\hline
Nonstandard Haar & 0.1400 &  & \\
\hline
Walsh Max & 0.005  & Walsh SVD &0.005\\
\hline
\end{tabular}
\end{center}
\end{table}
We present the $p$-values for different algorithms. If the p-value is less than $5\%$, 
one can confidently claim that the coordinates are not independent (rejecting the 
independence hypothesis). Clearly, some basis choices can better identify whether 
the data are independent or not. The Max-FiT algorithm using the nonstandard Haar wavelet 
basis (``Nonstandard Haar") is not competitive due to the local properties 
of the nonstandard Haar basis functions. 
Also, as the nonstandard Haar basis is not constructed using
the tensor product of one-dimensional basis functions, one cannot apply the SVD-FiT 
algorithm to ``Nonstandard Haar" to get global information and improve the 
performance. However, the ``sparse grid" nonstandard Haar wavelet basis may allow 
the study of high-dimensional data sets located on a low-dimensional
manifold. Therefore, an ongoing research topic is how to define a 
global statistic based on the adaptive nonstandard Haar wavelet basis 
for higher-dimensional data sets. Results along this direction
will be reported in future publications.

In our final example for the independence test, we study whether there exists a connection 
between the daily percentage alteration in temperature in New York City and that in the Dow
Jones Industrial Average (DJI) in 2019. 
\begin{figure}[htbp]
    \centering
    \includegraphics[width = 3in,height = 2.5in]{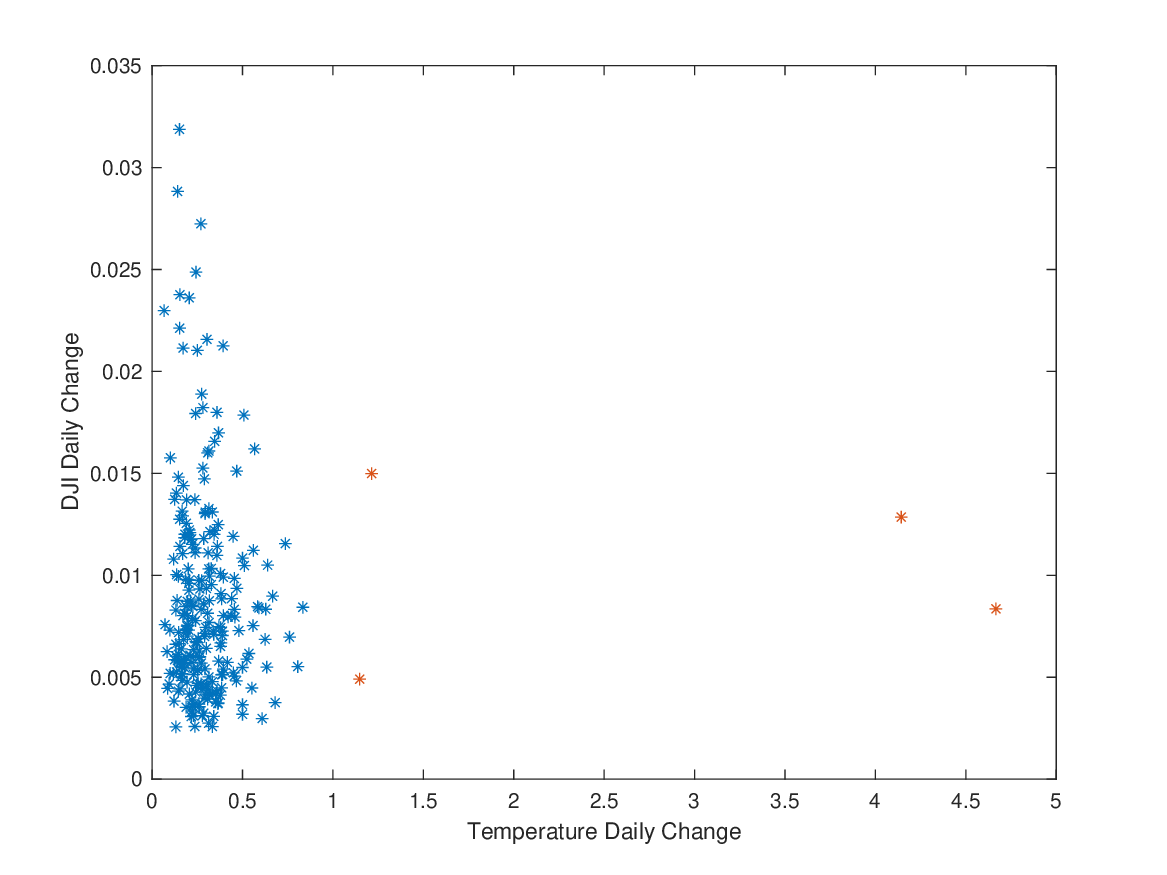}
    \caption{Temperature alternation vs DJI alternation. 
    \label{fig:temp_DJI}}
\end{figure}
For Day $i$, we define 
$x_i =$(highest temperature - lowest temperature)/(lowest temperature) and 
$y_i =$(highest DJI - lowest DJI )/(lowest DJI). In Table~\ref{tab:DJI}, 
we show the $p$-values from different algorithms. We observe significant 
differences in the results when different basis choices are used. The 
``smooth" and ``global" basis choices (Fourier series and Legendre polynomials) 
suggest that the data are independent, 
while methods using the Haar and Walsh wavelet bases seem to suggest that 
the data are related.
\begin{table}[htbp]
\caption{p-values for testing dependency of changes in temperature and in DJI index.}
	\label{tab:DJI}
\begin{center}
\begin{tabular}{|c|c|c|c|}
\hline
 & p-value &  & p-value\\
\hline
Fourier Max & 0.275 & Fourier SVD & 0.510\\
\hline
Haar Max & 0.060 & Haar SVD &0.015\\
\hline
 Legendre Max& 0.315  & Legendre SVD  &0.635\\
\hline
Nonstandard Haar & 0.005 &  & \\
\hline
Walsh Max & 0.060  & Walsh SVD & 0.015\\
\hline
\end{tabular}
\end{center}
\end{table}
One conjecture is that there exists a ``local" relation in the variables which states 
that when there is an extremely large percentage change in the temperature, the change
in the DJI tends to be milder. This local relation can be captured by the Walsh and 
Haar basis functions, but is smoothed out by the more global and smooth Legendre and 
Fourier bases. To validate this conjecture, we remove the 4 data points with the most 
significant temperature variations with $x_i>1$ (red points in Figure~\ref{fig:temp_DJI}), 
and re-calculate the p-values for the new data set. The results are shown in 
Table~\ref{tab:DJI_n}. All algorithms suggest that the new data sets are now more likely to be independent. 
\begin{table}[htbp]
\caption{p-values for testing dependency of changes in temperature and in DJI index after removing 4 points.}
	\label{tab:DJI_n}
\begin{center}
\begin{tabular}{|c|c|c|c|}
\hline
 & p-value &  & p-value\\
\hline
Fourier Max & 0.310 & Fourier SVD & 0.520\\
\hline
Haar Max & 0.205 & Haar SVD &0.110\\
\hline
 Legendre Max& 0.215  & Legendre SVD  &0.500\\
\hline
Nonstandard Haar & 0.100 &  & \\
\hline
Walsh Max & 0.285 & Walsh SVD & 0.110\\
\hline
\end{tabular}
\end{center}
\end{table}
This example demonstrates the challenges of accurate independence testing when the noise level is 
high (see also the differences in statistical power for large values $\sigma$
in Figure~\ref{fig:l_p}) or the relationships are complicated.

\subsection{Finding Low-Dimensional Structures}
\begin{figure}[htbp]
    \centering
    \includegraphics[width = 5in,height = 5in]{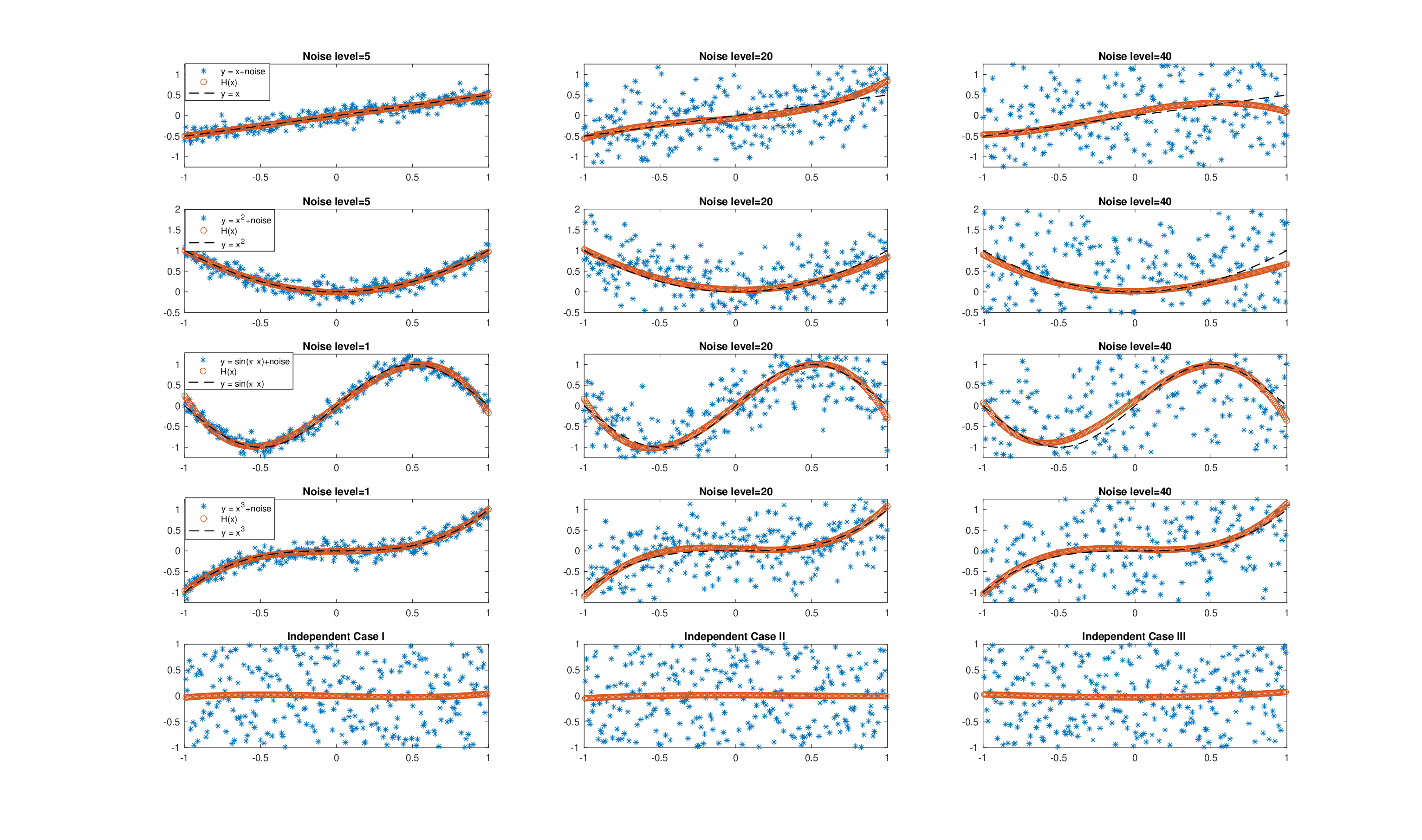}
    \caption{Bump Hunting results for different noise levels. \label{fig:b_h}}
\end{figure}

We first apply the bump hunting algorithm (Algorithm 3) to $5$ manufactured examples shown
in Figure~\ref{fig:b_h}. 
The (dependent) random data sets (blue asterisks) are generated using 
$y = x +\epsilon_1$ (linear), $y = x^2 + \epsilon_2$ (parabolic), 
$y = \sin(x) + \epsilon_3$ (sine) and $y = x^3+\epsilon_4$ (cubic), respectively, 
where $\epsilon_i \sim l*N(0,1)$, $i = 1,2,3,4$. For each example, we consider three 
different noise levels given by $l=5$, $l=20$, and $l=40$. 
For comparison reasons, we also consider 3 independent data sets sampled from the 
distribution $U([-1,1]^2)$ shown in the last row of Figure~\ref{fig:b_h}.
Each data set contains a total number of $256$ data points.
In the bump hunting algorithm implementation, we use the Legendre polynomials so that the 
identified low-dimensional structure is described by a polynomial function. 
For all cases, the frequency domain-based bump hunting algorithm (Algorithm 3) performs 
well for finding the well-defined bump functions, and the computed bump function 
(red circles) remains close to the ground truth (black dashed lines) even at
high noise levels. 

\begin{figure}[htbp]
    \centering
    \includegraphics[width = 5in,height = 4in]{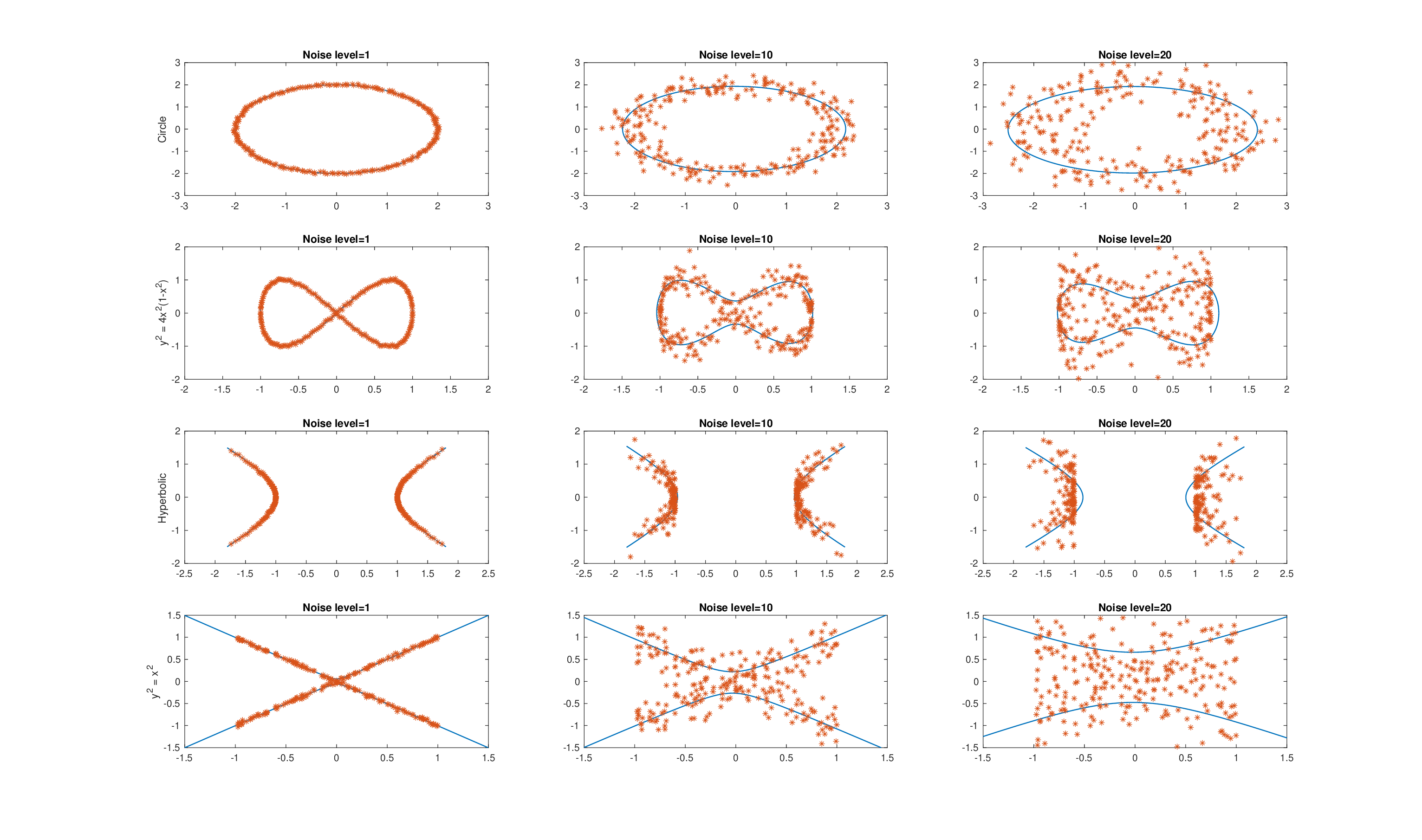}
    \caption{CCA-LDS results for data sets with different noise levels.
    \label{fig:cca1_3}}
\end{figure}

Next, we study the algebraic relations in the form of $H(x)=G(y)$ derived by solving
the specific optimization formulation discussed in Section~\ref{sec:bump}.
As polynomial approximations and the zeros of polynomial systems (algebraic varieties) 
are well-studied by the mathematics community, we therefore consider the manufactured 
data sets generated using polynomial relations listed in Table~\ref{tab:low-d}, 
with random noises $\epsilon_i \sim N(0,l/40)$, $i = 1,\cdots,6$ and 3 different 
noise levels $l=1$, $l=10$, and $l=20$, respectively. Each data set contains $256$ data 
points (red asterisks in Fig.~\ref{fig:cca1_3}).
\begin{table}[htbp]
\caption{Finding low-dimensional structures in the form of $H(x)=G(y)$.
	\label{tab:low-d} }
\begin{center}
\begin{tabular}{|c|c|c|c|}
\hline
 $H(x)=G(y)$ & x &y &$\theta$\\
\hline
$x^2+y^2=4$ & $x = 2\sin\theta + \epsilon_1$ & $y = 2\cos \theta+\epsilon_2$ &$[0,2\pi]$\\
\hline
 $y^2 = 4x^2(1-x^2)$& $x = \sin \theta $ & $y = \sin(2\theta)+\epsilon_3$& $ [0,2\pi]$\\
\hline
 $x^2 - y^2  = 1$& $x = \frac{\cos \theta}{\sqrt{\cos(2\theta)}}$ & $y 
%   = \frac{\sin \theta}{\sqrt{\cos(2\theta)}}+\epsilon_4$&$[-\frac{\pi}{4}+0.1,\frac{\pi}{4}-0.1] \cup  [\frac{3\pi}{4}+0.1,\frac{5\pi}{4}-0.1]$ \\
  = \frac{\sin \theta}{\sqrt{\cos(2\theta)}}+\epsilon_4$ & $[-0.69,0.68] \cup  [2.46,3.83]$ \\
\hline
$x^2 = y^2$ & $x\sim U(-1,1)$ & $y = \pm x+\epsilon_{5,6}$ & \\
\hline
\end{tabular}
\end{center}
\end{table}
We have applied both Algorithm 4.1 (CCA-LDS) and Algorithm 4.2 (SVD-LDS) to 
the cases listed in Table~\ref{tab:low-d}. As the results from Algorithm 4.1
and Algorithm 4.2 are very similar, we only plot the identified curves (blue curve) 
from Algorithm 4.1 in Fig.~\ref{fig:cca1_3}. For low noise levels, the algorithms
accurately capture the underlying low-dimensional manifolds. However, 
as the noise level increases, the accuracy of the identified low-dimensional curves
deteriorates.

\begin{figure}[htbp]
    \centering
    \includegraphics[width = 5in]{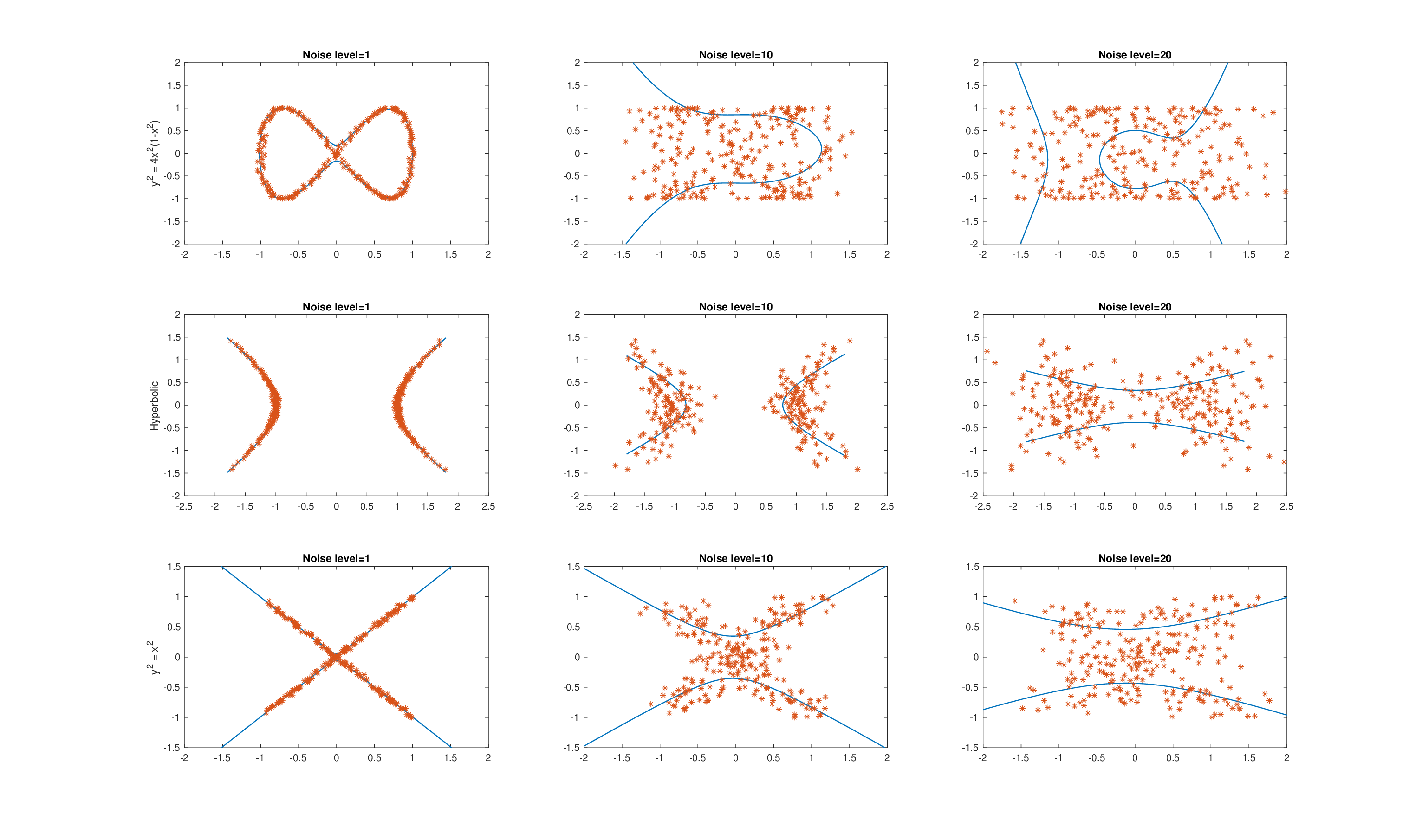}
    \caption{CCA-LDS results for data sets with noise only added to $x$.
    \label{fig:cca1_4}}
\end{figure}

Note that for the last three cases in Table~\ref{tab:low-d}, we only added noise to the $y$-variable. 
In the next experiment, we show how the results change when the noise is added only to 
the $x$-variable as demonstrated in Table~\ref{tab:low-d2}.
\begin{table}[htbp]
\caption{Adding noise only to the $x$-variable}
	\label{tab:low-d2}
\begin{center}
\begin{tabular}{|c|c|c|c|}
\hline
  $H(x)=G(y)$  & x &y &$\theta$\\
\hline
 $y^2 = 4x^2(1-x^2)$& $x = \sin \theta + \epsilon_1$ & $y = \sin(2\theta)$& $ [0,2\pi]$\\
\hline
 $x^2 - y^2  = 1$& $x = \frac{\cos \theta}{\sqrt{\cos(2\theta)}} + \epsilon_2$ & $y 
   = \frac{\sin \theta}{\sqrt{\cos(2\theta)}}$ & $[-0.69,0.68] \cup  [2.46,3.83]$ \\
\hline
$x^2 = y^2$ & $x= \pm y+\epsilon_{3,4}$ & $y\sim U(-1,1)$ & \\
\hline
\end{tabular}
\end{center}
\end{table}
We use the same random noise $\epsilon_i$ settings as those in Table~\ref{tab:low-d}.
The generated data sets (red asterisks) and identified low-dimensional structures (blue curves) 
are presented in Figure~\ref{fig:cca1_4}.
For low noise settings, Algorithms 4.1 and 4.2 can still correctly capture the ground truth. 
However, for high noise levels, as the relation $y^2 = 4x^2(1-x^2)$ is more nonlinear in $x$ 
than in $y$, the identified curves do not correctly capture 
the underlying curve used to generate the data sets.

\begin{figure}[htbp]
    \centering
    \includegraphics[width = 5in]{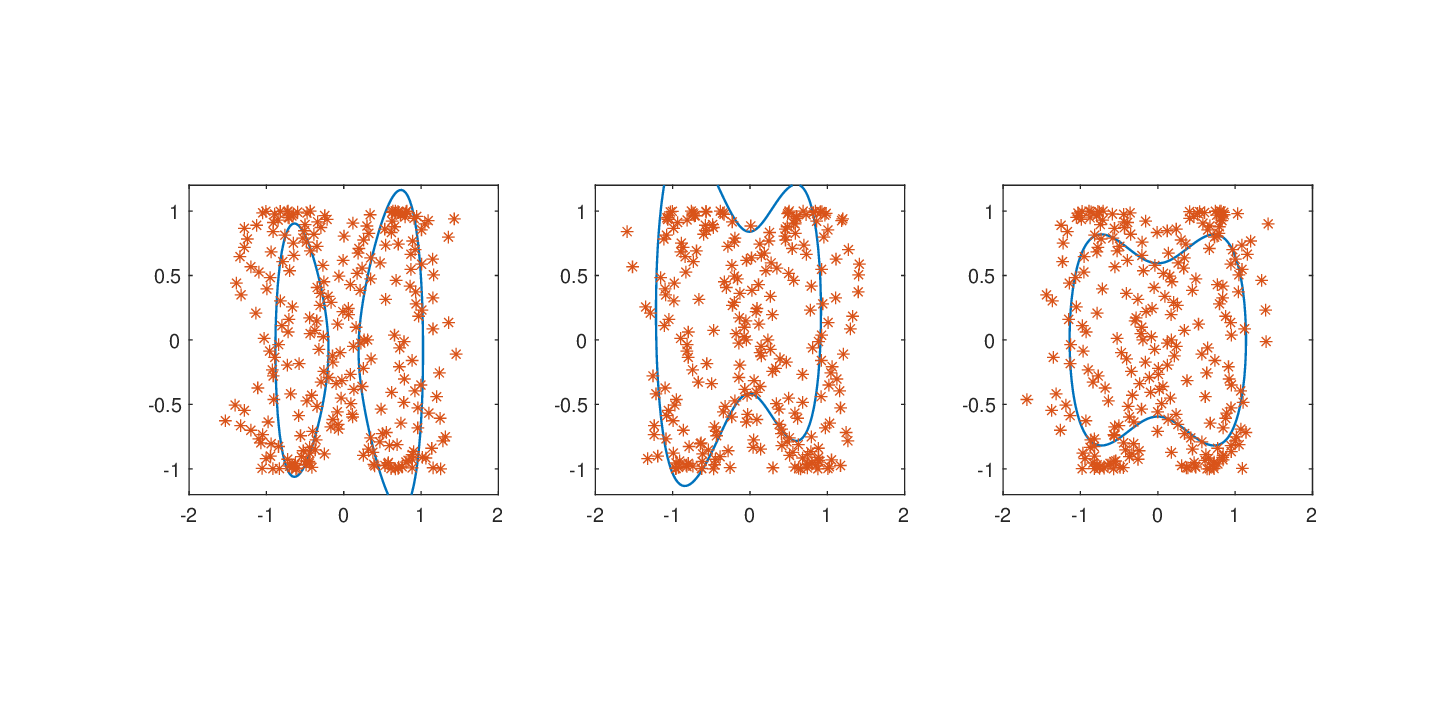}
    \caption{CCA-LDS results for relation $y^2 = 4x^2(1-x^2)$ with noise only added to the $x$-variable, 
	noise level $l=10$. Three data sets are obtained using different seeds for the random number generator.
    \label{fig:cca1_5}}
\end{figure}

We also find that the current Algorithms 4.1 and 4.2 are sensitive to random noises for cases
with high noise levels.
In Figure~\ref{fig:cca1_5}, we consider the relation $y^2 = 4x^2(1-x^2)$ with noise only 
added to the $x$-variable as presented in Table~\ref{tab:low-d2}. For the noise level $l=10$,
we generate three data sets using different seeds for the random number generator.
One can observe the differences in the identified low-dimensional manifolds (blue curves) although the same
relation and noise level are used to generate the random data sets. This also reveals the challenges
in numerically identifying the underlying nonlinear low-dimensional structures in high noise settings.

\begin{figure}[htbp]
    \centering
    \includegraphics[width = 4in]{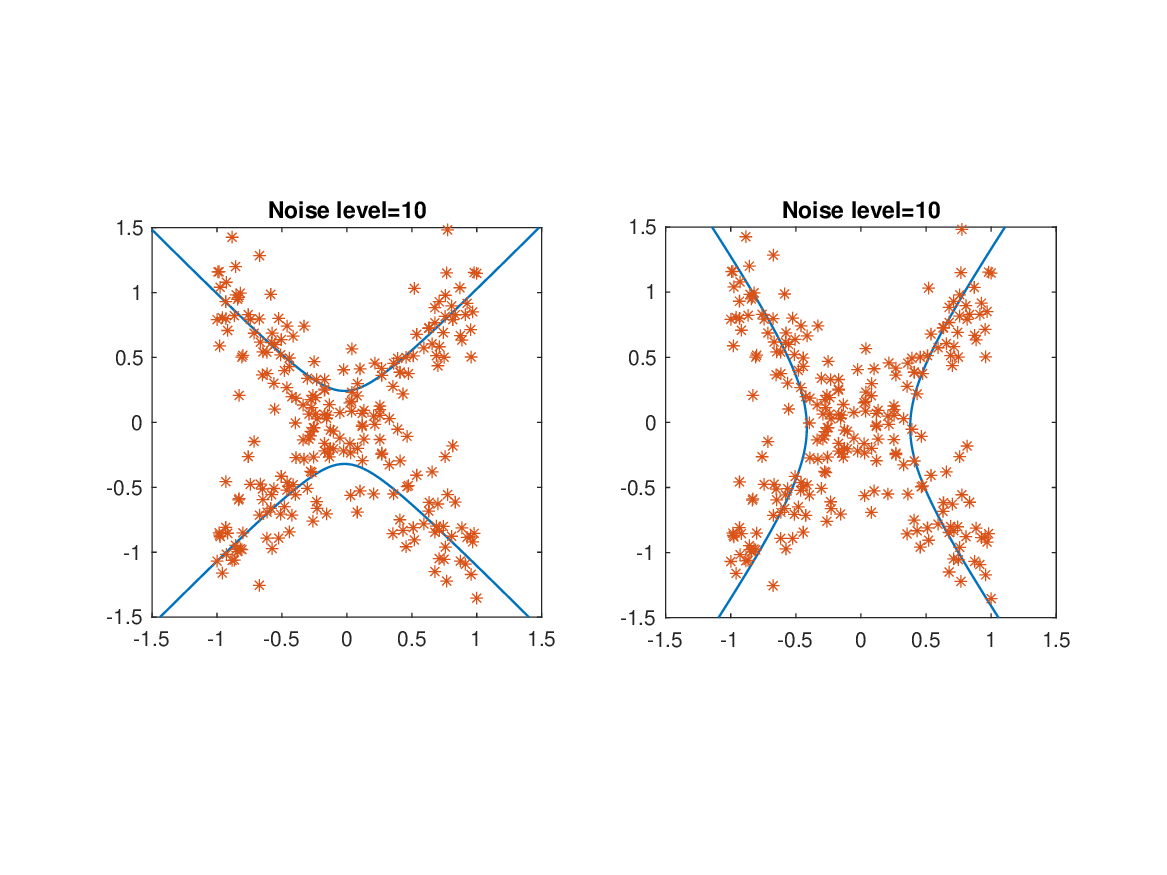}
    \caption{Identified low-dimensional structures using different linear fit strategies.
    \label{fig:cca1_6}}
\end{figure}
Finally, we show the impact of different linear fit strategies, $G(y) = \lambda H(x) + \gamma$ and 
$\tilde{\lambda} G(y) + \tilde{\gamma} = H(x)$, on the identified low-dimensional manifold in
Figure~\ref{fig:cca1_6}. We consider the symmetric relation $x^2 = y^2$ where the noise is added to
the $y$-variable only. The noise level is set to $l=10$. The data set contains $256$ points.
In the left panel, we show the identified curve when the standard linear fit $G(y) = \lambda H(x) + \gamma$
is applied. In the right panel, we show the result when the alternative linear fit 
$\tilde{\lambda} G(y) + \tilde{\gamma} = H(x)$ is used. 

%{\noindent \bf Comment:} 
Our preliminary research on the algebraic relations in the form of $H(x)=G(y)$ has revealed the challenges 
in finding low-dimensional structures in high-dimensional data sets. The current optimization formulation
has some nice existence, uniqueness, and regularity properties, and it can be solved efficiently using 
existing $L_2$ norm-based singular value decomposition or canonical-correlation analysis tools. 
However, the current formulation has stability issues and the identified low-dimensional 
algebraic structure appears to be sensitive to noise structures 
(e.g., how the noise is added and noise levels). Some of these issues are being 
studied for more general algebraic relations $H(x,y)=0$ and results will be 
reported in the future.

%% file: sec05Summary0503.tex
\section{Summary and Future Work}
\label{sec:summary}
In this paper, by introducing different orthonormal basis functions 
and applying fast transformation algorithms, we construct the frequency 
domain expansion coefficients and higher-order moments of the probability density 
function of physical domain data. We show how the new set of frequency domain
data can be utilized to address some fundamental questions arising in data analysis, 
including testing the dependence between variables in the data set and 
identifying any low-dimensional algebraic relations in dependent data.
Numerical experiments are presented to demonstrate the performance of the resulting 
numerical algorithms.
 
Revisiting the fundamental data analysis questions from the frequency/spectral domains 
and applied analysis perspectives has provided us with powerful tools to address 
some interesting and challenging questions. One of our current projects is to generalize
the frequency domain statistical analysis techniques to high-dimensional data sets, by
properly choosing an appropriate set of adaptive basis functions for high-dimensional
data with low-rank or low-dimensional compressible features. Another project is to 
design improved statistics that better distinguish the null distribution from the 
alternative distribution. Results along these directions will be reported in the 
future.

% One of our current research projects is to design better strategies to rank the coefficients of the frequency domain expansion in the dependency tests, which may be useful to help identify a global ``skeleton" of the frequency domain data. Another project is to using a combination of different basis, or more general frames (redundant bases) to further improve the statistical power of existing frequency domain methods. In many engineering applications, it is common to use both function values and their derivative values to more effectively analyze a particular function. Results along these directions will be reported in the future.

\section*{Acknowledgments}
We thankfully acknowledge the support from the National Science Foundation under 
Grants \#2152289 (HB, JH, JS, KZ, WZ), \#2012451 (SC,JH,JS), and \#2152070 (YW). 
Any opinions, findings, and conclusions or recommendations expressed in this 
paper are those of the authors and do not necessarily reflect the views of the 
National Science Foundation.